\documentclass[12pt]{amsart}
\usepackage{enumerate,amstext,amsfonts,amssymb,amscd,amsbsy,amsmath}
\usepackage{ifthen,multicol,extarrows,fullpage,anyfontsize,mathtools,multirow}
\usepackage{parskip}
\usepackage{lscape}
\usepackage{enumitem}
\usepackage{color,tikz}
\usepackage{amsthm}
\usepackage{latexsym}
\usepackage[all]{xy}
\usepackage[bookmarks, colorlinks=true, linkcolor=blue, citecolor=blue, urlcolor=blue]{hyperref}
\usepackage[backrefs,lite,alphabetic]{amsrefs} 
\usepackage{tikz-cd}
\usepackage{pgfplots}

\usepackage{heatMapMacros}

\usepackage{eucal}

\numberwithin{equation}{section}
\newtheorem{lemma}[equation]{Lemma}

\newtheorem{conj}[equation]{Conjecture}

\newtheorem{claim*}{Claim}

\newtheorem{defn}[equation]{Definition}

\newtheorem{question}[equation]{Question}

\theoremstyle{remark}
\newtheorem{example}[equation]{Example}
\newtheorem{remark}[equation]{Remark}

\newcommand{\Proj}{\operatorname{Proj}}

\newcommand{\Tor}{\operatorname{Tor}}

\newcommand{\Sym}{\operatorname{Sym}} 
\newcommand{\GL}{\mathbf{GL}}

\newcommand{\Alt}{\bigwedge\nolimits}

\newcommand{\bd}{\mathrm D}

\newcommand{\PP}{\mathbb P}

\newcommand{\QQ}{\mathbb Q}
\newcommand{\ZZ}{\mathbb Z}
\newcommand{\KK}{\mathbb K}

\newcommand{\bdelta}{\boldsymbol{\delta}}
\newcommand{\CC}{\mathbb C}
\newcommand{\ba}{\mathbf a}
\newcommand{\bb}{\mathbf b}

\renewcommand{\bd}{\mathbf d}
\newcommand{\be}{\mathbf e}

\newcommand{\bzero}{\mathbf 0}

\newcommand{\bS}{\mathbf{S}}

\newcommand{\cO}{\mathcal O}

\newcommand{\zero}{\mathbf 0}

\usepackage{rotating}

\title{Syzygies of $\PP^{1}\times \PP^{1}$: data and conjectures}

\author{Juliette Bruce}
\address{Department of Mathematics, University of California, Berkeley, CA}
\email{\href{mailto:juliette.bruce@berkeley.edu}{juliette.bruce@berkeley.edu}}
\urladdr{\url{https://juliettebruce.github.io}}

\author{Daniel Corey}
\address{Department of Mathematics, University of Wisconsin, Madison, WI}
\email{\href{dcorey@math.wisc.edu}{dcorey@math.wisc.edu}}
\urladdr{\url{https://sites.google.com/site/dcorey2814/home}}

\author{Daniel Erman}
\address{Department of Mathematics, University of Wisconsin, Madison, WI}
\email{\href{mailto:derman@math.wisc.edu}{derman@math.wisc.edu}}
\urladdr{\url{http://math.wisc.edu/~derman/}}

\author{Steve Goldstein}
\address{Botany Department and Department of Biostatistics and Medical Informatics, University of Wisconsin, Madison, WI}
\email{\href{mailto:sgoldstein@wisc.edu}{sgoldstein@wisc.edu}}

\author{Robert P. Laudone}
\address{Department of Mathematics, University of Wisconsin, Madison, WI}
\email{\href{laudone@wisc.edu}{laudone@wisc.edu}}
\urladdr{\url{https://www.math.wisc.edu/~laudone/}}

\author{Jay Yang}
\address{School of Mathematics, University of Minnesota, Minneapolis, MN}
\email{\href{mailto:jkyang@umn.edu}{jkyang@umn.edu}}
\urladdr{\url{http://www-users.math.umn.edu/~jkyang/}}
\date{\today}


\begin{document}
\thanks{
JB was partially supported by the National Science Foundation under Award Nos. DMS-1502553, DMS-1440140, and NSF MSPRF DMS-2002239. DC received support from NSF-RTG grant 1502553. DE received support from NSF grant DMS-1601619 and DMS-1902123.  RL received support from NSF grant DMS-2001992 and DMS-1502553. JY received support from NSF grants DMS-1502553 and DMS-1745638.
}

\begin{abstract}
We provide a number of new conjectures and questions concerning the syzygies of $\PP^1\times \PP^1$. The conjectures are based on computing the graded Betti tables and related data for large number of different embeddings of $\PP^1\times \PP^1$. These computations utilize linear algebra over finite fields and high-performance computing.
 \end{abstract}

\maketitle

\section{Introduction}
While syzygies are a much-studied topic in algebraic geometry and commutative algebra, the Betti tables for varieties of dimension $\geq 2$ remain largely mysterious.  For instance, the Betti table of $\PP^2$ under the $d$-uple Veronese embedding is only fully understood for $d\leq 6$  \cite{bruceErmanGoldsteinYang18, wcdl}, and there is not yet even a conjectural picture for the values of such Betti tables.  One obstacle to developing such a conjecture is a lack of data: for the $d$-uple embedding of $\PP^2$, the required number of variables grows like $d^2$, and so free resolution computations tend to overflow memory.

In~\cite{bruceErmanGoldsteinYang18},  the computation of syzygies was approached via an alternate method.  Instead of using symbolic Gr\"obner basis methods to compute a minimal free resolution, we computed the Betti numbers via the cohomology of the Koszul complex.  In essence, this swapped a symbolic computation for a massive linear algebra computation.  (See \S\ref{sec:background} for the theoretical background on this approach.)  This reduced the computation to a number of individual rank computations, one for each multigraded Betti number, and then we performed those computations using high-throughput computations.

The present work has three foci: we improve the framework for this alternate approach to Betti numbers; we apply it to the case of $\PP^1\times \PP^1$ to generate a wealth of new data; and we use that data to offer new conjectures and questions about the syzygies of $\PP^1\times \PP^1$. 

\subsection{Overview of the computation}

For any $\bd=(d_1,d_2)\in \ZZ_{>0}^2$, we can embed $\iota_{\bd}\colon \PP^1\times \PP^1\to \PP^{(d_1+1)(d_2+1)-1}$ by  the complete linear series for $\cO_{\PP^1\times \PP^1}(\bd)$, and we want to understand the syzygies of this image. Following a philosophy implicit in Green's foundational work on syzygies~\cite{green-II}, and echoed in later results on asymptotic syzygies~\cite{einLazarsfeld12,eel-quick}, we will study the syzygies of not only the structure sheaf, but also of the pushforward of various line bundles.  In particular, our goal is to compute the syzygies of $\iota_{\bd *}\cO_{\PP^1\times \PP^1}(\bb)$ for as many choices of $\bd$ and $\bb$ as possible.

Depending on the grading group or equivariant structure under consideration, we can represent these Betti numbers in a multitude of ways.  See \S\ref{sec:background} for a summary of notation. 

Our main computation involves the $\ZZ^4$-graded Betti numbers. There are $\approx d_1d_2$ entries of the Betti table which could be nonzero, and each of those entries will involve at most $\approx d_1^3d_2^3$ distinct $\ZZ^4$ multidegrees.  However, by using known vanishing and duality results, accounting for symmetry, and applying elementary results on the relationship between Betti numbers and Hilbert function, we can shrink down to a much smaller number of matrices, which we refer to as the {\bf relevant range}, and which are sufficient to determine all of the Betti numbers.  (See \S\ref{sec:main computation} for details on the relevant range.)

The main computation involves computing the ranks of all of the matrices from the relevant range. The rank of each matrix can be computed in parallel, allowing us to leverage high throughput computational resources. In addition, some of the matrices are quite massive, and we thus require huge amounts of memory for those particular matrices. 

\begin{figure}
	\includegraphics[scale=0.6]{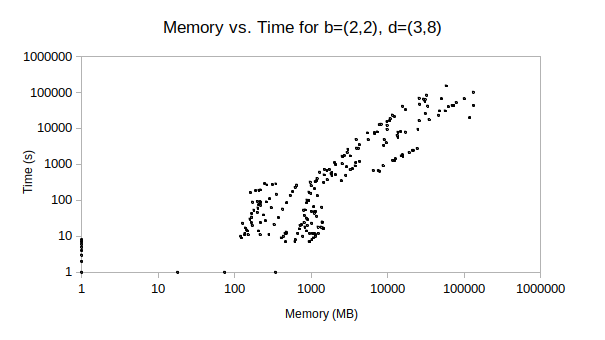}
	\caption{Memory vs. time to compute ranks of matrices for $\bb = (2,2)$, $\bd = (3,8)$} 
	\label{fig:memoryTimeScatter}
\end{figure}

For concreteness, let us consider our largest complete computation, which is the case $\bb=(2,2)$ and $\bd = (3,8)$.  The relevant range involves 1130 matrices, the largest of which is $2,124,896 \times 3,719,448$, and Figure~\ref{fig:memoryTimeScatter} provides a scatterplot of the time and memory involved in computing the ranks of those matrices.  Only a handful of cases took over a day.

\subsection{Computational improvements}
Our current work improves on the method of~\cite{bruceErmanGoldsteinYang18} in a number of ways.  Most notably, \cite{bruceErmanGoldsteinYang18} relied on floating-point rank computations of sparse real matrices, using a MATLAB implementation of the LU-algorithm; by contrast, our current work simply performs the computations over finite fields in MAGMA.  MAGMA recently introduced major improvements in their linear algebra of finite fields~\cite{steel}, which seemed to make these rank computations much faster than our previous method; see Figure~\ref{fig:LU vs MAGMA}.  

\begin{figure}
\[
\begin{tabular}{c | c | c}
Method&Average time per job (secs) & Max time (secs)\\ \hline
MatLab LU-algorithm over $\mathbb R$& 220&4735\\
MAGMA rank algorithm over $\mathbb F_{32003}$ & 7&99
\end{tabular}
\]
\caption{We compared floating-point LU-algorithm computations in MatLab with rank computations in MAGMA over the finite field $\mathbb F_{32003}$, for all of the multigraded matrices related to one individual Betti number.  This anecdotally suggests that MAGMA computations over finite fields are significantly faster, though we did not do any comprehensive testing.}
\label{fig:LU vs MAGMA}
\end{figure}

Moreover, this switch to working over finite fields enabled us to use exact calculations, eliminating the need for floating-point approximations.  While an exact computation over a finite field will not necessarily agree with the exact computation over $\mathbb Q$, there are only finitely many primes where the computations could disagree, and these discrepancies seem to rarely arise for reasonably large primes.  This switch to working over finite fields thus had a significant downstream effect: the main computations in ~\cite{bruceErmanGoldsteinYang18} introduced some numerical errors as $\bd$ grew larger, requiring the use of representation theoretic techniques to detect these errors.  By contrast, our finite field computations produced no such numerical errors, and we were able to produce Schur functor decompositions without the need for the sort of ``error correction'' from~\cite[\S5]{bruceErmanGoldsteinYang18}.

\subsection{New Data}
After computing the multigraded ranks for the relevant range, we process the data into usable formats.  The rank computations quickly yield $\ZZ^4$-multigraded Betti numbers, but most mathematical conjectures focus on either the standard $\ZZ$-graded Betti numbers or on the underlying $\GL_2\times \GL_2$-Schur modules.  We convert into those formats and encode all of the results into a Macaulay2 package for ease of use.

In total, we compute complete Betti tables for just shy of 200 total pairs of $\bb$ and $\bd$. See \S\ref{sec:data} and Table~\ref{table:computedData} for more details on the data.

\subsection{Conjectures}
Based on the data we computed, we develop a number of new conjectures, and we provided evidence in support of some previous conjectures.

\begin{figure}
\includegraphics[scale=0.7]{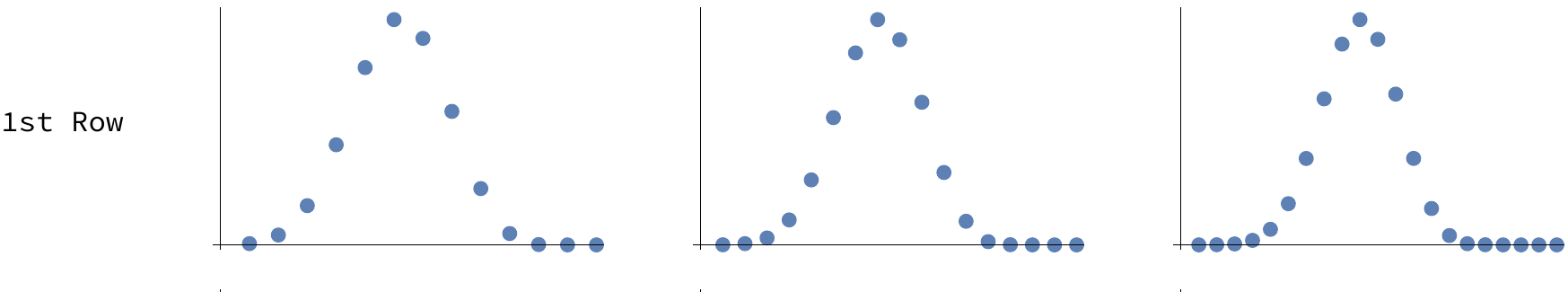}
\caption{Here we plot the Betti numbers of the first row of the Betti tables of $\PP^1\times \PP^1$ embedded by $(d,3)$, for $d=3,4$ and $5$.  They appear to resemble a normal distribution, as predicted by a conjecture of~\cite{EEL}.}
\label{fig:3d normal distribution}
\end{figure}

We first examine the quantitative behavior of the standard grade Betti numbers, with conjectures in \S\ref{sec:qual conjectures} that address unimodality properties of the Betti numbers and various statistics.  In addition, we consider our data in relation to a conjecture of Ein, Erman, and Lazarsfeld that, for large values of $\bd$, the Betti numbers in any given row of the Betti table should behave like a binomial distribution~\cite[Conjecture~B]{EEL}.  A theorem of Bruce~\cite[Theorem~A]{bruce-semiample} implies that the first row of the Betti table for $\PP^1\times \PP^1$ and line bundles $(2,d_2)$ satisfies exactly this behavior as $d_2\to \infty$.  Our data provides further support for the normal distribution behavior suggested by the conjecture, and seems to show this behavior even for the low values of $d_2$ for which we have data.  See Figure~\ref{fig:3d normal distribution}.

In \S\ref{sec:rep theory conjectures}, we consider several conjectures related to the $\GL_2\times\GL_2$ structure of these syzygies.  This includes precise conjectures on the Schur functor decomposition of certain entries; an analysis of the shapes of partitions that arise; and a discussion of ``redundant'' representations.  

In \S\ref{sec:boij sod}, we present a collection of conjectures involving the Boij-S\"oderberg decompositions of these Betti tables. In particular, we provide a complete conjectural description of the Boij-S\"oderberg coefficients of the homogeneous coordinate ring of $\PP^1\times \PP^1$ embedded by $\cO_{\PP^1\times \PP^1}(2,d_2)$

\section*{Acknowledgments} 
We thank Erika Pirnes for her contributions to early versions of some of our code.
We thank UW-Madison Math Department, Center for High Throughput Computing, John Canoon, Claudiu Raicu, Greg Smith, and Allan Steel. The first author is grateful for the support of the Mathematical Sciences Research Institute in Berkeley, California, where she was in residence for the 2020-2021 academic year.  The computer algebra systems Magma and Macaulay2 provided valuable assistance throughout our work~\cites{M2,MAGMA}.

\section{Background and Notation}\label{sec:background}
Throughout this section, we work over an arbitrary field $\KK$.   Our convention will be to write integer vectors using boldface, as in $\bd\in\ZZ^2$, and to specify the coordinates as $\bd=(d_1,d_2)$. We let $\bzero=(0,0)\in \ZZ^{2}$. 

As we are interested in the syzygies of $\PP^1\times \PP^1$ throughout we let $S=\KK[x_{0},x_{1},y_{0},y_{1}]$ be the corresponding polynomial of over a field $\KK$. When viewed as the Cox ring of $\PP^1\times \PP^1$~\cite{cox95}, the ring $S$ inherits a $\ZZ^2$-bigrading given by $\deg x_0 = \deg x_{1} = (1,0)\in \ZZ^2$ and $\deg y_{0} = \deg y_{1} = (0,1) \in \ZZ^{2}$.  The ring $S$ also admits a $\ZZ^4$-multigrading given by setting the degree of each variable to be a generator of $\ZZ^4$, e.g. $\deg(x_0) = (1,0,0,0)$ and $\deg(x_1)=(0,1,0,0)$ and so on.

\subsection{Standard graded Betti numbers}\label{subsec:betti}
The syzygies of $\PP^1\times \PP^1$ under various embeddings come from studying Segre-Veronese modules of $S$.  Given $\bd\in \ZZ^2_{>0}$ and $\bb\in \ZZ^{2}$ the Segre-Veronese module is
\begin{equation*}
S(\bb;\bd) \coloneqq \bigoplus_{k\in \ZZ} S_{k\bd+\bb}.
\end{equation*}
Since $k\bd+\bb$ determines a ray in $\ZZ^2$ as $k$ varies in $\ZZ$, $S(\bb;\bd)$ is naturally a $\ZZ$-graded module over the polynomial $R=\Sym S_{\bd}$. When $\bb=\zero$ the module $S(\zero;\bd)$ is isomorphic to the homogeneous coordinate ring of $\PP^1\times\PP^1$ embedded by $\cO_{\PP^1\times\PP^1}(\bd)$ into the projective space $\PP^{(\bd_1+1)(\bd_2+1)-1}=\Proj R$.  If $\bb\ne \zero$, then $S(\bb;\bd)$ is naturally isomorphic to the section module of a pushforward of a line bundle; specifically, $S(\bb;\bd)$ is the $R$-module associated to the sheaf $(\iota_{\bd})_{*} \cO_{\PP^1\times\PP^1}(\bb)$.  As noted in the introduction, while our primary interest is in the syzygies of the homogeneous coordinate rings $S(\zero;\bd)$, past work shows that studying the syzygies of other line bundles is often helpful in providing a more uniform picture \cite{einLazarsfeld12,eel-quick,green-II}.

The Betti numbers of a graded $R$-module $M$ are defined as $\beta_{i,j}(M) = \dim_{\KK} \Tor^R_i(M,\KK)_j$, which denotes the degree $j$ part of the $\Tor_i$-module.  For convenience, when studying the Betti numbers $S(\bb;\bd)$, we will omit reference to the ambient polynomial ring $R$, and write $\beta_{i,j}(S(\bb;\bd)) = \beta_{i,j}(\PP^1\times\PP^1,\bb;\bd)$.  The Betti numbers of a graded module are often computed using a minimal free resolution~\cite{eisenbud-syzygies,M2}.  However, an alternate characterization of the Betti numbers, via Koszul cohomology, is more relevant for our computational approach.    

The Koszul complex of $S(\bb;\bd)$ over the ring $R$ is the complex:
\begin{equation*}
\begin{tikzcd}[column sep = 2.5em]
\cdots \rar{}&\Alt^{1}R^{(\bb_1+1)(\bb_2+1)}\otimes S(\bb;\bd)\rar{}&\Alt^{0}R^{(\bb_1+1)(\bb_2+1)}\otimes S(\bb;\bd),
\end{tikzcd}
\end{equation*}
which is naturally $\ZZ$-graded since since $R$ is $\ZZ$-graded.  Given a pair of integers $(p,q)$, we can analyze the cohomology of the degree $p+q$ strand of this complex, in homological degree $p$.  This will be denoted by $K_{p,q}(\PP^1\times\PP^1,\bb;\bd)$.\footnote{We remark that $K_{p,q}$ and $\beta_{i,j}$ provide two different notations for similar invariants, though $K_{p,q}$ is a vector space whereas $\beta_{i,j}$ is an integer; both are commonly used in the literature.  We will primarily use the $K_{p,q}$-notation, however the conversion between the two notations is given by the simple rule $\dim K_{p,q}\leftrightarrow \beta_{p,p+q}$.}  It can be computed explicitly as the middle cohomology of the following complex:
\begin{equation}\label{eq:kozul-S}
\begin{tikzcd}[column sep = 2.5em]
\cdots \rar{}& \Alt^{p+1}S_{\bd}\otimes S_{(q-1)\bd+\bb}\rar{\partial_{p+1,q-1}}&\Alt^{p}S_{\bd}\otimes S_{q\bd+\bb}\rar{\partial_{p,q}}&\Alt^{p-1}S_{\bd}\otimes S_{(q+1)\bd+\bb}\rar{}&\cdots
\end{tikzcd}
\end{equation}
where the differentials are given by
\begin{align*}
\partial_{p+1,q-1}\left(m_0\wedge m_1\wedge\cdots\wedge m_p\otimes f\right)&=\sum_{i=0}^{p}(-1)^i m_0\wedge m_1\wedge\cdots\wedge \hat{m}_i\wedge \cdots\wedge m_p\otimes m_if \\
\partial_{p,q}\left(m_1\wedge m_2\wedge\cdots\wedge m_p\otimes f\right)&=\sum_{i=1}^{p}(-1)^i m_1\wedge m_2\wedge\cdots\wedge \hat{m}_i\wedge \cdots\wedge m_p\otimes m_if.
\end{align*}
In other words, instead of computing all of the Betti numbers simultaneously via a minimal free resolution, we can compute each Betti number individually using the complex of vector spaces in \eqref{eq:kozul-S}.  This, in essence, turns a problem of symbolic algebra into a (massive but largely distributable) problem in linear algebra.

\subsection{Multigraded Betti numbers}\label{subsec:multi}
By incorporating the $\ZZ^4$-grading on $S$, we can subdivide the problem even further and obtain the $\ZZ^4$-graded Betti numbers.  For a multidegree $\mathbf{e}\in \ZZ^4$, we define $\beta_{i,\mathbf{e}}(\PP^1\times\PP^1,\bb;\bd)=\beta_{i,\mathbf{e}}(S(\bb;\bd)) = \dim_{\KK}\Tor^R_i(S(\bb;\bd),\KK)_\mathbf{e}$.  This is well defined because both $R$ and $S(\bb;\bd)$ inherit $\ZZ^4$-multigradings from $S$.  From the Koszul cohomology perspective, the Koszul complex of $S(\bb;\bd)$ over $R$ is also homogeneous with respect to the $\ZZ^4$-grading.  Thus, we can analyze the cohomology of the degree $\mathbf{e}$-strand, which provides our method for computing $\beta_{i,\mathbf{e}}(\PP^1\times\PP^1,\bb;\bd)$.

\subsection{Schur functor decomposition}\label{subsec:schur}
The action of $\GL_2\times \GL_2$ on $\PP^1\times \PP^1$ turns the vector space $K_{p,q}(\PP^1\times\PP^1,\bb;\bd)$ into a $\GL_2\times \GL_2$-representation.  We can therefore decompose $K_{p,q}(\PP^1\times\PP^1,\bb;\bd)$ into a direct sum of irreducible $\GL_2\times \GL_2$-representations.  These irreducible representations have the form $\bS_\lambda \otimes \bS_\mu$, where $\lambda,\mu$ are partitions with length $\leq 2$. See \cite[Exercise 2.36]{FultonHarris} for background. For brevity, we write $\bS_{(a,b,c,d)}$ for the Schur module $\bS_{(a,b)}\otimes \bS_{(c,d)}$.

\begin{example}
	Let $\bb = (2,2)$ and $\bd = (3,3)$. The Betti table for $S(\bb;\bd)$ is
{\footnotesize
	\begin{align*}
	\begin{array}{ccccccccccccccc}
	&0&1&2&3&4&5&6&7&8&9&10&11&12&13 \\
	0:&9&108&585&1872&3861&5148&4026&1080&\mathbf{165}&\cdot&\cdot&\cdot&\cdot&\cdot \\
	1:&\cdot&\cdot&\cdot&\cdot&\cdot&165&1080&4026&5148&3861&1872&585&108&9 
	\end{array}
	\end{align*}
}
	The bold entry in the Betti table tells us that $\dim K_{8,0}(\PP^1\times \PP^1,\bb;\bd) = \beta_{8,8}(\PP^1\times \PP^1,\bb;\bd) = 165$. Viewed as $\GL_2\times\GL_2$-representation, $K_{8,0}(\PP^1\times \PP^1,\bb;\bd)$ decomposes as
	\begin{equation*}
	K_{8,0}(\PP^1\times \PP^1, \bb;\bd) \cong
	\bS_{(17, 9, 17, 9)} \oplus \bS_{(16, 10, 16, 10)} \oplus \bS_{(15, 11, 15, 11)} \oplus \bS_{(14, 12, 14, 12)}\oplus \bS_{(13, 13, 13, 13)}.
	\end{equation*}
The dimensions of these Schur modules are $81, 49, 25, 9$ and $1$, respectively.
\end{example}

\subsection{Koszul Duality}
Using duality of Koszul cohomology groups (see, for instance~\cite[Duality Theorem (2.c.9)]{green-I}), we can derive data for more values of $\bb$ and $\bd$, as we now explain. Given $\bb$ we define its \textit{Koszul dual} as $\bb' := \bd - \bb - (2,2)$.  We have
\begin{equation*}
K_{p,q}(\PP^1\times \PP^1,\bb';\bd) \cong K_{(d_1+1)(d_2+1)-3 -p, 2-q}(\PP^1\times \PP^1,\bb;\bd)
\end{equation*}
as vector spaces. Visually, this means that the Betti table for $(\bb';\bd)$ is obtained by rotating the Betti table for $(\bb;\bd)$ by $180^{\circ}$. We will illustrate this phenomenon in Example \ref{ex:koszulDual}.  Note that $(d_1+1)(d_2+1)-3$ is the codimension of $\PP^1\times \PP^1$ in the embedding by $\bd$.
The duality also applies to the Schur functor decomposition via the following formula.   To phrase this, we need some more notation.  Let
\[
\alpha:=\left(\tfrac{(d_1+1)(d_2+1)d_1 -2}{2}, \tfrac{(d_1+1)(d_2+1)d_1 -2}{2}, \tfrac{(d_1+1)(d_2+1)d_2 -2}{2}, \tfrac{(d_1+1)(d_2+1)d_2 -2}{2}\right).
\]
Given any $w = (w_0,w_1,w_2,w_3)\in \ZZ^4$ we write  $w^{\text{opp}}=(w_1,w_0,w_3,w_2)$ and we choose $w'$ so that $w+(w')^{\text{opp}} = \alpha$.   The multiplicity of the Schur functor $S_w$  in $K_{p,q}(\bb;\bd)$ equals the multiplicity of the Schur functor $S_{w'}$ in the dual Koszul cohomology group $K_{(d_1+1)(d_2+1)-3 -p, 2-q}(\PP^1\times \PP^1,\bb';\bd)$, where $\bb'$ is defined as above.

\begin{example}
	\label{ex:koszulDual}
	Let $\bb = (0,0)$ and $\bd = (3,3)$. The Betti table for $S(\bb;\bd)$ is
	{\footnotesize
	\begin{align*}
	\begin{array}{ccccccccccccccc}
	&0&1&2&3&4&5&6&7&8&9&10&11&12&13\\
	0:&1&\cdot&\cdot&\cdot&\cdot&\cdot&\cdot&\cdot&\cdot&\cdot&\cdot&\cdot&\cdot&\cdot\\
	1:&\cdot&87&676&2691&6864&12155&15444&14157&9152&3861&780&\mathbf{22}&\cdot&\cdot\\
	2:&\cdot&\cdot&\cdot&\cdot&\cdot&\cdot&\cdot&\cdot&\cdot&\cdot&165&144&39&4
	\end{array}
	\end{align*}
    }
	\noindent The bold entry in the Betti table tells us that $\dim K_{11,1}(\PP^1\times \PP^1,\bb;\bd) = \beta_{11,12}(\PP^1\times \PP^1,\bb;\bd) = 22$. Viewed as $\GL_2\times\GL_2$-representation, $K_{11,1}(\PP^1\times \PP^1,\bb;\bd)$ decomposes as
	\begin{equation*}
	K_{11,1}(\PP^1\times \PP^1,\bb;\bd) \cong
	\bS_{(23, 13, 18, 18)} \oplus \bS_{(18, 18, 23, 13)}.
	\end{equation*}
	
	The Koszul dual pair to $(\bb;\bd)$ is  $\bb' = (1,1)$ and $\bd' = (3,3)$. The Betti table for $S(\bb';\bd')$ is
	{\footnotesize
	\begin{align*}
	\begin{array}{ccccccccccccccc}
	&0&1&2&3&4&5&6&7&8&9&10&11&12&13\\
	0:&4&39&144&165&\cdot&\cdot&\cdot&\cdot&\cdot&\cdot&\cdot&\cdot&\cdot&\cdot\\
	1:&\cdot&\cdot&\mathbf{22}&780&3861&9152&14157&15444&12155&6864&2691&676&87&\cdot\\
	2:&\cdot&\cdot&\cdot&\cdot&\cdot&\cdot&\cdot&\cdot&\cdot&\cdot&\cdot&\cdot&\cdot&1
	\end{array}
	\end{align*}
    }
	\noindent We see that this Betti table is exactly that corresponding to $S(\bb;\bd)$ rotated by $180^{\circ}$. The $(11,1)$ entry for $S(\bb;\bd)$ corresponds to $(2,1)$ for $S(\bb';\bd')$. Viewed as $\GL_2\times\GL_2$-representation, $K_{2,1}(\PP^1\times \PP^1,\bb;\bd)$ decomposes as
	\begin{equation*}
	K_{2,1}(\PP^1\times \PP^1,\bb;\bd) \cong
	\bS_{(10, 0, 5, 5)} \oplus \bS_{(5, 5, 10, 0)}
	\end{equation*}	
\end{example}

\section{Computed Data}\label{sec:data}
Using the algorithms outlined in Section~\ref{sec:main computation} we computed the Betti tables, $\ZZ^4$-multigraded Betti numbers, and Schur functor decompositions for over 150 distinct pairs  $(\bb;\bd)$, including $27$ distinct $\bd$-values. In Table \ref{table:computedData}, we list, for each $\bd$, the number of $\bb$'s for which we have complete data.  For comparison: ~\cite{bruceErmanGoldsteinYang18} computed similar data for $\PP^2$ for about $15$ distinct pairs $(\bb;\bd)$, which included $5$ distinct $\bb$ values; and ~\cite{wcdl}, which only considered the case $\bb=\zero$, computed data for $\PP^2$ for $5$ distinct $\bd$ values.   There appears to be no significant computational work on syzygies for $\PP^1\times \PP^1$, although \cite{lemmensP1P1} does construct a non-minimal resolution.  In other words, these computations represent a significant contribution to the available syzygy data for $\PP^1\times \PP^1$ specifically, as well as for toric surfaces more generally.

\begin{table}[tbh!]
	\begin{tabular}{|ll|lllllllll|}
		\hline
		&   &   &   &   &  & $d_2$ &   &   &   &    \\
		&   &2&3&4&5&6&7&8&9&10 \\ \hline
		& 2 &3&6&8&10&12&14&13&6&6   \\
		$d_1$ & 3 &$\cdot$&6&12&15&13&12&8&4&2    \\
		& 4 &$\cdot$&$\cdot$&9&14&9&5&1&1&0    \\
		& 5 &$\cdot$&$\cdot$&$\cdot$&1&1&1&1&0&0   \\
		\hline
	\end{tabular}
\caption{For each $\bd$, the number of $\bb$ for which we compute the Betti tables, $\ZZ^4$-multigraded Betti numbers, and Schur decompositions}
\label{table:computedData}
\end{table}

\begin{remark}
	In Table \ref{table:computedData}, for the symmetric cases $\bd = (d,d)$, we only record $\bb=(b_1,b_2)$ with $b_1\leq b_2$ for which we have data. For example, when $\bd=(2,2)$, we only count the cases $\bb = (0,0)$, $(0,1)$, and $(1,1)$; we do not include $(1,0)$.
\end{remark}

\section{Main Computation}\label{sec:main computation}

Broadly speaking, our approach to computing the Betti table, $\ZZ^4$-multigraded Betti numbers, and Schur functor decompositions for a given pair $(\bb;\bd)$ proceeds as follows:
\begin{enumerate}
	\item {\bf Reduction to the relevant range:}  By combining a computation of the multigraded Hilbert series with known  vanishing results for syzygies (relying primarily on Castelnuovo-Mumford regularity), we conclude that a small subset of the Betti numbers determines all of the Betti numbers. This smaller subset is the relevant range, and is the focus of our computations.
	\item {\bf Constructing the matrices in the relevant range:}  We follow the ideas in~\cite{bruceErmanGoldsteinYang18} to efficiently construct and store the matrices from the relevant range.
	\item  {\bf High throughput rank computations:}  We use distributed high throughput computation to find the ranks of all the matrices in the relevant range.  These computations are done via linear algebra over the finite field $\mathbb F_{32003}$ in MAGMA.  This is by far the most computationally intensive aspect.
	\item  {\bf Post-processing:}    Using standard ideas from representation theory, we convert the multigraded Betti number into Schur functor decompositions.
\end{enumerate}

While the techniques here are broadly similar to those in \cite{bruceErmanGoldsteinYang18}, which focused on computing syzygies of Veronese embeddings of $\PP^{2}$, the passage from $\PP^2$ to $\PP^1\times \PP^1$ requires new code in each step and we further refine this implementation and approach. The most significant distinction is in the third step abvoe:  the core algorithm in the current work uses linear algebra over finite fields, whereas in ~\cite{bruceErmanGoldsteinYang18} it used floating-point computations.

\subsection{Relevant Range}

We expedite our computations significantly by utilizing the fact that for many values of $p$ and many multidegrees $\ba$, the multigraded Betti number $\beta_{p,\ba}(\PP^1\times\PP^1,\bb;\bd)$ is determined entirely by the $\ZZ^4$-multigraded Hilbert series of $S(\bb;\bd)$.  In the following lemma, we use vector notation $t^{\ba}:=t_0^{a_0}t_1^{a_1}t_2^{a_2}t_3^{a_3}$ if $\ba=(a_0,a_1,a_2,a_3)$.

\begin{lemma}
The $\ZZ^4$-multigraded Hilbert series of $S(\bb;\bd)$ is a rational function of the form:
$
A(t_0,t_1,t_2,t_3) / B(t_0,t_1,t_2,t_3)
$ where
\[
A = \sum_{p,\ba} \beta_{p,\ba}(\PP^1\times \PP^1, \bb;\bd)t^{\ba} \text{ and } B = \prod_{\bb \in \mathbb N^4, b_0+b_1=d_1, b_2+b_3=d_2} (1-t^{\bb})
\]
\end{lemma}
The proof is nearly identical to that of~\cite[Lemma 3.1]{bruceErmanGoldsteinYang18}, so we omit it.

With this in mind, our main computations reduce to determining the ranks $\partial_{p,q}$ for $p,q$ in what we call the relevant range. 

\begin{defn}
Fixing $\bb$ and $\bd$ we define the relevant range to be the set of pairs $(p,q)$ such that $K_{p,q}(\PP^1\times\PP^1,\bb;\bd)\neq0$ and either $K_{p-1,q+1}(\PP^1\times\PP^1,\bb;\bd)\neq0$  or $K_{p+1,q-1}(\PP^1\times\PP^1,\bb;\bd)\neq0$.
\end{defn}

In general we determine the relevant range by finding the smallest $p$ such that $K_{p,q}(\PP^{1}\times\PP^{1},\bb;\bd)\neq0$ and then applying duality (see \cite[Proposition 3.5]{einLazarsfeld12}). When $\bb=\bzero$ the only case of interest is $q=1$, and we find the smallest $p$ such that $K_{p,1}(\PP^{1}\times\PP^{1},\bzero;\bd)\neq0$ via  \cite[Theorem~1.4]{wcdl}. When $\bb \neq \bzero$ we determine the relevant range using the fairly coarse vanishing bounds from~\cite[Proposition 5.1]{einLazarsfeld12}.  While a sharper bound on the relevant range would allow us to compute ranks for many fewer matrices, we found that in practice, these potentially extraneous matrices did not cause any bottlenecks in the actual computation.

An algorithm entirely analogous to~\cite[Algorithm 3.3]{bruceErmanGoldsteinYang18} enables us to efficiently compute the multigraded Betti numbers outside of the relevant range.

\subsection{Constructing the matrices in the relevant range}

 After computing the relevant range and the relevant multidegrees, this data is fed to the code to compute the matrices representing the differentials in the relevant range. We first use the $\mathfrak{S}_2\times \mathfrak{S}_2$-symmetries of the multidegrees to restrict to those multidegrees $(a,b,c,d)$ where $a\geq b$ and $c\geq d$. As in \cite{bruceErmanGoldsteinYang18} we use duality for Koszul cohomology groups to reduce the number of matrices we compute~\cite[Theorem~2.c.6]{green-I}. 
 Unfortunately unlike in the case of the Veronese, the bi-graded structure means that it is not possible to use this duality to reduce to a finite set of non-redundant Betti tables.

When constructing the matrices, we use the fact that all of the maps $(\partial_{p,q})_{\ba}$ correspond to submatrices of the boundary map $d_p:\bigwedge^p S_{\bd}\rightarrow \bigwedge^{p-1} S_{\bd}$. In particular, $(\partial_{p,q})_{\ba}$ is given by restricting to the submatrix $d_{p,\leq \ba}$ given by those entries in degrees $\leq \ba$. However, instead of storing the map $d_p$ we simply use this fact to compute all of the various $(\partial_{p,q})_{\ba}$ for all multidegrees at once. This was implemented as it was found that as the degrees got larger, more of the entries in the $d_p$ matrix correspond to multidegrees that are not in the relevant range. This is entirely analogous to~\cite[\S4.1]{bruceErmanGoldsteinYang18}, which provides further details.  In Appendix \ref{appendix:matrices}, we list the number of matrices we must compute and the largest such matrix.

\begin{example} 
  For $\bd=(3,8)$, $\bb=(2,2)$, the full computation of which is discussed in more detail in Example~\ref{ex:b22d38}, it took a modern laptop computer, 5min 25sec to compute all the relevant matrices, entailing a total of 1130 matrices, taking a total of 13GB of space. The single largest matrix had 16,999,168 non-zero entries.
\end{example}

\subsection{High Throughput Computations}
The rank computations can be efficiently distributed over numerous different computers. We implemented these computations using high throughput computing via HTCondor on the  University of Wisconsin--Madison Mathematics department computer servers.
Many of the matrices are small, and hence do not require much memory to compute the rank. Because our hardware grid has fewer nodes with large amounts of available RAM, the initial submissions are allocated a small amount of RAM (e.g. 2GB). For the jobs that fail, we resubmit with a larger memory allocation, and repeat this  process until the computation terminates.

\begin{figure}[htbp]
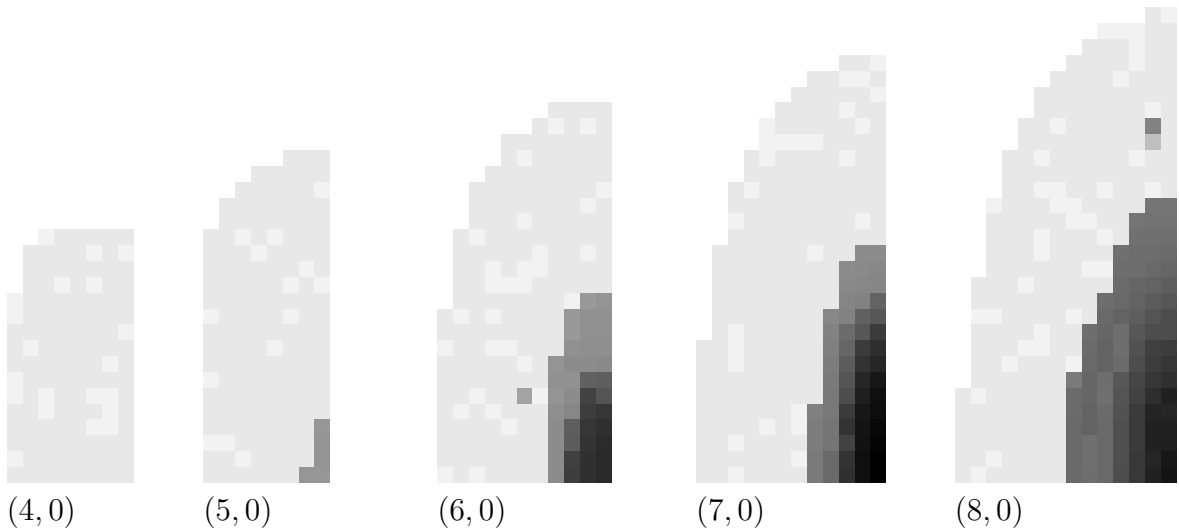

	\centering
	\begin{minipage}[b]{0.15\textwidth}
	\begin{stackColor}[271]
0	0	15	25	25	25	25	25\\
0	25	25	25	25	15	25	15\\
0	25	25	25	25	25	25	25\\
0	25	25	15	25	15	25	25\\
15	25	25	25	25	25	25	25\\
15	25	25	25	25	25	25	25\\
25	25	25	25	25	25	25	15\\
25	15	25	25	25	25	25	25\\
25	25	25	25	25	25	15	25\\
15	25	25	25	25	25	25	25\\
15	25	15	25	25	15	15	25\\
25	25	15	25	25	25	15	25\\
25	25	25	25	25	15	15	25\\
25	25	25	25	25	25	25	25\\
15	25	25	25	25	25	25	25\\
25	25	25	25	25	25	25	25
	\end{stackColor} \\  {\centering $(4,0)$}
\end{minipage}
\begin{minipage}[b]{0.18\textwidth}
		\begin{stackColor}[271]
	0	0	0	0	0	0	0	0\\
	0	0	0	0	0	25	25	25\\
	0	0	0	25	25	25	25	25\\
	0	0	25	25	25	25	25	15\\
	0	25	25	25	25	25	25	25\\
	0	25	25	25	25	25	25	25\\
	25	25	15	25	15	25	25	25\\
	25	25	25	15	25	25	25	25\\
	25	25	25	25	25	25	15	25\\
	25	25	25	25	25	15	25	15\\
	25	25	25	25	25	25	25	25\\
	15	25	25	25	25	15	25	25\\
	25	25	25	25	25	25	25	25\\
	25	25	25	25	15	25	25	25\\
	25	25	25	25	25	25	25	25\\
	15	25	25	25	25	25	25	25\\
	25	25	25	25	25	25	25	25\\
	25	25	25	25	25	25	25	25\\
	25	25	25	25	25	25	25	113\\
	15	15	25	25	25	25	25	114\\
	25	25	15	25	25	25	25	115\\
	25	25	25	25	25	25	113	113
\end{stackColor} \\$(5,0)$
	\end{minipage}
	\begin{minipage}[b]{0.2\textwidth}
		\begin{stackColor}[271]
	0	0	0	0	0	0	0	0	0	0	0\\
	0	0	0	0	0	0	0	0	0	0	0\\
	0	0	0	0	0	0	0	25	25	25	25\\
	0	0	0	0	0	0	25	15	25	15	25\\
	0	0	0	0	25	25	25	25	25	25	25\\
	0	0	0	0	25	15	25	25	25	25	25\\
	0	0	0	25	25	25	25	25	25	25	25\\
	0	0	25	25	25	25	25	25	25	25	15\\
	0	0	25	25	25	25	25	25	25	15	25\\
	0	0	25	25	25	15	25	25	25	25	25\\
	0	25	15	25	25	25	25	25	25	25	25\\
	0	25	25	25	25	25	15	25	25	15	25\\
	0	25	25	15	25	15	15	25	25	25	25\\
	0	25	25	15	15	15	25	25	25	25	25\\
	0	25	25	25	25	25	25	25	15	110	114\\
	25	15	25	15	25	25	25	25	110	117	118\\
	25	25	25	25	25	25	25	25	111	119	118\\
	15	25	25	15	15	25	25	25	118	121	120\\
	25	25	25	25	25	15	25	115	124	122	126\\
	25	25	25	25	25	25	25	118	117	170	172\\
	25	25	15	25	25	99	15	116	124	201	183\\
	25	15	25	15	25	25	25	121	128	208	216\\
	25	25	25	25	15	25	25	122	173	215	222\\
	25	25	25	25	25	25	25	122	175	217	226\\
	25	25	25	25	25	25	25	123	205	222	229\\
	15	25	15	25	25	25	25	123	206	221	229
\end{stackColor} \\$(6,0)$
	\end{minipage}
	\begin{minipage}[b]{0.2\textwidth}
		\begin{stackColor}[271]
					0	0	0	0	0	0	0	0	0	0	0	0\\
		0	0	0	0	0	0	0	0	0	0	0	0\\
		0	0	0	0	0	0	0	0	0	0	0	0\\
		0	0	0	0	0	0	0	0	0	25	25	15\\
		0	0	0	0	0	0	0	25	25	15	15	25\\
		0	0	0	0	0	0	25	25	25	25	25	15\\
		0	0	0	0	0	25	25	25	25	15	25	25\\
		0	0	0	0	15	25	25	25	25	25	15	25\\
		0	0	0	0	15	15	15	15	25	25	25	25\\
		0	0	0	25	15	25	25	25	25	15	25	25\\
		0	0	0	25	25	25	25	25	25	25	25	25\\
		0	0	25	15	25	25	25	25	25	25	25	25\\
		0	0	25	25	25	25	25	25	25	25	25	25\\
		0	0	15	25	25	25	25	25	25	25	15	25\\
		0	0	25	25	25	25	25	25	25	25	25	25\\
		0	25	25	25	25	25	25	15	25	25	122	124\\
		0	25	25	25	25	25	25	25	25	123	126	128\\
		0	25	25	25	25	25	25	25	25	127	130	132\\
		0	25	25	25	25	25	25	25	25	131	134	168\\
		0	25	25	25	25	25	25	25	129	134	169	181\\
		0	25	15	25	25	25	25	25	133	158	180	210\\
		25	25	15	25	25	25	25	25	136	164	209	224\\
		25	25	15	25	25	25	25	25	139	175	218	233\\
		25	25	25	25	25	25	25	25	142	183	232	248\\
		25	25	25	25	25	25	25	25	144	190	240	257\\
		25	25	25	25	25	25	15	138	146	209	247	259\\
		25	25	25	25	15	25	15	139	147	219	251	264\\
		25	25	15	25	25	25	25	140	149	202	253	268\\
		25	25	25	25	25	25	15	141	161	219	256	271\\
		25	25	15	25	25	25	25	141	162	220	257	271
		\end{stackColor} \\$(7,0)$
	\end{minipage}
	\begin{minipage}[b]{0.2\textwidth}
		\begin{stackColor}[271]
			0	0	0	0	0	0	0	0	0	0	0	0	0	0\\
			0	0	0	0	0	0	0	0	0	0	0	0	0	0\\
			0	0	0	0	0	0	0	0	0	0	0	0	0	0\\
			0	0	0	0	0	0	0	0	0	0	0	0	0	0\\
			0	0	0	0	0	0	0	0	0	0	0	0	25	15\\
			0	0	0	0	0	0	0	0	0	15	15	15	25	25\\
			0	0	0	0	0	0	0	0	25	25	25	15	25	25\\
			0	0	0	0	0	0	0	25	25	25	25	15	25	25\\
			0	0	0	0	0	0	25	25	25	25	15	25	25	25\\
			0	0	0	0	0	25	25	25	25	25	25	25	25	25\\
			0	0	0	0	0	25	25	25	25	25	25	25	15	25\\
			0	0	0	0	25	25	15	25	25	25	25	25	134	25\\
			0	0	0	0	25	25	25	25	25	25	25	25	68	25\\
			0	0	0	25	25	15	25	25	25	25	25	15	25	25\\
			0	0	0	25	25	25	25	25	25	25	25	25	25	25\\
			0	0	0	25	25	15	15	25	25	15	25	25	15	25\\
			0	0	15	25	25	25	15	15	25	25	25	25	147	148\\
			0	0	25	25	25	25	25	15	15	25	25	149	150	151\\
			0	0	25	25	25	15	25	25	25	15	25	152	154	154\\
			0	0	25	25	25	15	25	25	25	25	153	155	157	158\\
			0	0	25	25	25	25	25	25	15	15	156	155	161	163\\
			0	25	25	25	15	25	25	25	25	25	159	163	167	171\\
			0	25	25	25	25	25	25	25	25	157	162	167	176	180\\
			0	15	15	25	25	25	25	15	25	160	166	173	184	190\\
			0	25	25	25	25	15	25	25	157	163	156	180	192	190\\
			0	25	25	25	15	15	25	25	159	166	157	186	204	206\\
			0	25	25	15	25	25	25	15	161	166	158	191	207	215\\
			0	25	25	25	25	25	25	141	163	156	179	198	216	222\\
			25	15	25	25	25	25	25	158	165	153	182	202	224	237\\
			25	25	25	25	25	25	25	159	167	159	185	211	234	232\\
			25	25	25	25	25	25	25	160	167	161	187	211	237	238\\
			25	25	25	25	25	25	25	161	154	161	189	211	236	239\\
			25	25	15	25	25	25	25	161	154	162	189	213	233	249\\
			25	15	25	25	25	25	25	161	156	160	192	213	238	252
		\end{stackColor} \\$(8,0)$
	\end{minipage}
		\caption{Memory usage to compute ranks for each multidegree for $\bb = (2,2)$, $\bd = (3,8)$.  Plots are arranged left to right $(p,q) = (4,0)$ through $(8,0)$.} \label{fig:heatMapMultidegree}
\end{figure}

\begin{example}
	\label{ex:b22d38}
	In this example, we provide  a detailed analysis of how we determine the Betti table for $\bd=(3,8)$ and $\bb=(2,2)$, one of our larger computations. There are only two rows $q=0,1$, and 34 columns; we display the first several columns below.
	
	{\footnotesize
		\begin{align*}
		\begin{array}{ccccccccccc}
		&0&1&2&3&4&5&6&7&8&9\\
		0:&9&258&3465&28512&156546&568620&1210506&697680&203490&\cdot\\
		1:&\cdot&\cdot&\cdot&1050&28476&498498&5444400&41855958&194378184&671067540\\
		\end{array} \cdots		
		\end{align*} }
	
	\noindent The relevant range is $(p,0)$ for $4\leq p \leq 8$ and $(p,1)$ for $3\leq p \leq 7$. Because $K_{p,0} - K_{p-1,1} $ is determined by the Hilbert function of the module, we need only compute one of $K_{p,0}$ or $K_{p-1,1}$, and we compute the former.  To that end, we form the matrices $(\partial_{p,0})_{\ba}$ and $(\partial_{p+1,-1})_{\ba}$ for $4\leq p \leq 8$ and compute their ranks. Fortunately, $(\partial_{p+1,-1})_{\ba}=0$.  After accounting for  $\mathfrak{S}_2\times \mathfrak{S}_2$-symmetry, we are left to compute ranks of  $1130$ matrices, the largest of which is $2,\!124,\!896 \times 3,\!719,\!448$.
    In this case, up to symmetry there were 39788 multidegrees with non-zero entries in the Betti table. For these entries, in absence of the consideration about relevant ranges, to compute these entries would have required the computation of at least $81,\!437$ matrices.

The amount of RAM and time used in the rank calculation is recorded in Figure \ref{fig:memoryTimeScatter}.  The vast majority of matrices require less than 1MB of RAM and 10 seconds. Figure \ref{fig:memAndTimeVp} has two plots displaying the average and maximum memory, resp. time, needed to compute the ranks of the matrices $(\partial_{p,0})_{\ba}$ as a function of $p$.

\begin{figure}[h]
	\centering
	\begin{minipage}[b]{0.45\textwidth}
	\includegraphics[scale=0.55]{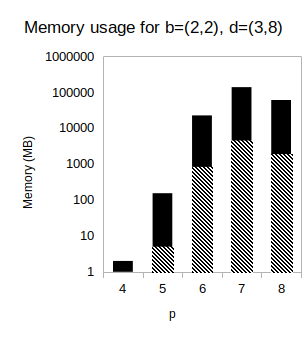}
	\end{minipage}
	\begin{minipage}[b]{0.45\textwidth}
		\includegraphics[scale=0.55]{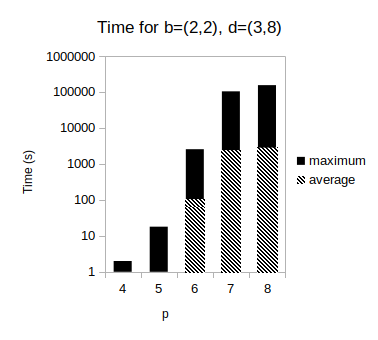}
	\end{minipage}
	\caption{Memory and time to compute ranks of matrices for $\bb = (2,2)$, $\bd = (3,8)$ and $q=0$} \label{fig:memAndTimeVp}
\end{figure}

Figure \ref{fig:heatMapMultidegree} illustrates how memory usage varies with multidegree for each $(p,0)$.  The plots are arranged left to right $(p,q) = (4,0)$ through $(8,0)$. Here is how to interpret these plots. Within each plot, each square represents a multidegree, and its color measures the memory usage:  light gray is 0 GB and black reaches the maximum of 132 GB of RAM. Because of the $\mathfrak{S}_2\times \mathfrak{S}_2$-symmetry, we need only consider the multidegrees $(a,b,c,d)$ satisfying $a+b=26$, $a\leq b$ and $c+d=66$, $c\leq d$. Each row has $(a,b)$ constant, each column has $(c,d)$ constant, and $a$, resp. $c$, increases in the downward, resp. left, direction.
\end{example}

Here are some take-aways from this example. We see that the amount of memory and time needed to compute ranks of matrices comprising the differential $\partial_{p,q}$ grows as $p$ moves towards the center of the Betti table. Nevertheless, for a fixed $(p,q)$, nearly all of the matrices $(\partial_{p,q})_{\ba}$ require minimal memory and time. The $(\partial_{p,q})_{\ba}$ that require the most resources are those for which $\ba=(a,b,c,d)$ are balanced, i.e., for which $|a-b|$ and $|c-d|$ are minimized.

\begin{remark}
The fact that the most computationally intensive are those for which $\ba$ is balanced could allow one to potentially dig deeper into conjectures related to Schur functors.  Namely, the highest weight of a given Schur module tends to be quite unbalanced.  Given the parallel nature of these computations, one could potentially rule out the presence of certain Schur modules for many values of $\bb$ and $\bd$ for which a full computation would be impossible.
\end{remark}

\subsection{Post-processing}
Having computed all of the multigraded Betti numbers, we can easily combine the values to obtain the standard graded Betti numbers.  Obtaining the Schur functor decompositions is a bit more involved, though it is nearly identical to the process in \cite[\S 5.1]{bruceErmanGoldsteinYang18}. The main idea is once again a highest weight greedy algorithm. In the $\PP^2$ case, the authors were considering the decomposition as a $\GL_3$-module. In our case, we are considering the decomposition as a $\GL_2 \times \GL_2$-module. The irreducible polynomial representations of $\GL_2 \times \GL_2$ are products $\bS_\lambda \otimes \bS_\mu$ of Schur functors where $\lambda,\mu$ are partitions with length $\leq 2$. For further details, see \cite[Chapter 6, Exercise 2.36]{FultonHarris}. 

We order bi-partitions using the standard Lex order on $\ZZ^4$. That is, for two bi-partitions $(\lambda,\mu), (\nu,\eta)$ we say that $(\lambda,\mu) \leq (\nu,\eta)$ if $(\lambda_1,\lambda_2,\mu_1,\mu_2) \leq (\nu_1,\nu_2,\eta_1,\eta_2)$ in the standard Lex order on $\ZZ^4$. This gives us a well order on bi-partitions. In particular, we can select a largest element.

To decompose $K_{p,q}(\PP^1 \times \PP^1; \bb; \bd)$ into Schur functors we apply the Algorithm for Schur Functor Decomposition (see below).  
The algorithm terminates due to the semi-simplicity of $\GL_2 \times \GL_2$.  More specifically, semi-implicitly implies that there is a finite decomposition $K_{p,q}(\PP^1 \times \PP^1; \bb; \bd) \cong \bigoplus_{\lambda,\mu} (\bS_\lambda(\CC^2) \otimes \bS_\mu(\CC^2))^{\oplus c_{\lambda,\mu}}$ for some constants $c_{\lambda,\mu}$ only finitely many of which are nonzero. This means the multigraded Hilbert series $H$ in the algorithm above is a sum of Hilbert series corresponding to $(\bS_\lambda(\CC^2) \otimes \bS_\mu(\CC^2))$ which is just the product of the Hilbert series for $\bS_\lambda$ in the variables $t_0,t_1$ and the Hilbert series for $\bS_\mu$ in a second set of variables $t_2,t_3$.  The weight of the lex-leading monomial of the Hilbert series will always be a bi-partition, i.e. $\lambda_1 \geq \lambda_2$ and $\mu_1 \geq \mu_2$; and that monomial will correspond to the highest weight of some Schur modules appearing in the decomposition.  Thus, the algorithm uses the lex-leading monomial of the Hilbert series to iteratively pick off summands in the decomposition of $K_{p,q}(\PP^1 \times \PP^1; \bb; \bd)$.


\begin{figure}
\noindent {\bf Algorithm for Schur Functor Decomposition}
\begin{align*}
&Input: &&\beta_{p,\ba}(\bb;\bd) \; \text{for fixed} \; \bb;\bd,p \; \text{and all} \; \ba \in \ZZ_{\geq 0}^4 \; \text{with} \; |\ba| = (p+q)(d_1+d_2) + (b_1+b_2)\\
&Output: &&\text{A list $K$ of bi-partitions appearing in the Schur module decomposition} \\
&&&\text{of $K_{p,q}(\PP^1 \times \PP^1; \bb; \bd)$, with multiplicity.}\\
&Steps: &&L \coloneqq \{\ba \; | \; |\ba| = (p+q)(d_1+d_2) + (b_1+b_2)\} \; \text{and} \; H = \sum_{\ba \in L} \beta_{p,\ba}(\bb; \bd) \cdot t^{\ba}\\
&&&K = \{\}\\
&&& \text{While the coefficient of ${\rm lex}(H) > 0$ do:}\\
&&& \qquad \text{Let $(\lambda,\mu) = (\lambda_1,\lambda_2,\mu_1,\mu_2)$ be the weight of the lex-leading monomial in $H$}\\
&&& \qquad \text{Let $K = K \cup \{(\lambda,\mu)\}$}\\
&&& \qquad \text{Let $H$ equal $H$ minus the multigraded Hilbert series of $\bS_{\lambda}(\CC^2) \otimes \bS_{\mu}(\CC^2)$}.\\
&&& \text{Return $K$}.
\end{align*}
\label{alg:decomp}
\end{figure}

 \section{Qualitative Aspects of the Computed Data}\label{sec:qual conjectures}

\subsection{Unimodality}

Our data strongly suggests that several statistics associated with the syzygies of $\PP^{1}\times\PP^{1}$ are unimodal. More specifically, our data leads to the following conjecture. 

\begin{conj}\label{conj:unimodal}
For any $\bb$, if either $d_1$ or $d_2$ is sufficiently large, then each of the following functions is unimodal:
\begin{enumerate}
	\item  The standard graded Betti numbers in a single row: $i\mapsto \beta_{i,i+k}(\PP^1\times \PP^1,\bb;\bd)$ for any fixed $k$.\footnote{In the range of $\bb$ we have considered in this paper, these functions are only interesting for $k=0,1$ or $2$.}
	\item  The number of Schur functors with multiplicity appearing in a given row: fix some $q$ and consider $p\mapsto$ the total number of Schur functors, counted with multiplicity, appearing in $K_{p,q}(\PP^1\times \PP^1,\bb;\bd)$.
	\item  The largest multiplicity of Schur functors appearing in a given row:  fix some $q$ and consider $p\mapsto$ the largest multiplicity of a Schur functors appearing in $K_{p,q}(\PP^1\times \PP^1,\bb;\bd)$.

\end{enumerate}
\end{conj}

\begin{remark}
Our data also suggests that even the multigraded Betti numbers exhibit unimodality in certain ways, although in the multigraded setting there is no canonical choice for what one might expect to be unimodal. For example, fixing any multidegree $\be$ and our data suggests that $i\mapsto  \beta_{i,i\be}(\PP^1\times \PP^1,\bb;\bd)$. It would be interesting to explore other ways in which the multigraded Betti numbers might satisfy some sort of unimodality or concavity properties. Given the large number of possible multidegrees, such questions can be somewhat complex.
\end{remark}

Patterns similar to Conjecture~\ref{conj:unimodal} were observed for the Veronese syzygies of $\PP^{2}$  in \cite{bruceErmanGoldsteinYang18}*{Section 6.4}. Interestingly in this setting the authors observed that the function $p\mapsto$ the number of distinct Schur functors appearing in $K_{p,q}(\PP^2,\cO_{\PP^{2}}(b);\cO_{\PP^{2}}(d))$ appears to be unimodal (see \cite{bruceErmanGoldsteinYang18}*{Question~6.11.(2)}). By contrast, our data provides a large number of counterexamples to that for $\PP^{1}\times \PP^{1}$. More specifically, out of the rough 90 pairs of $\bb$ and $\bd$ that we tested, the number of distinct Schur functors appearing was not unimodal.

\begin{example}
Letting $\bb=\bzero$, $\bd=(3,4)$, and considering $q=1$ we see that the number of distinct Schur functors appearing in the decompositions of $K_{p,1}(\PP^1\times\PP^1,\bb;\bd)$ is
\[
(9, 26, 42, 52, 67, 71, 82, {\bf 80}, 87, {\bf 78}, 79, 63, 49, 5, 1),
\]
which is not unimodal. We see a similar failure of the number of distinct Schur functors appearing in the decompositions of $K_{p,1}(\PP^1\times\PP^1,\bb;\bd)$ when $\bb=\bzero$ and $\bd=(3,5)$: 
\[
(11, 32, 56, 67, 96, 101, 127, {\bf 125}, 146, {\bf 137}, 154, {\bf 135}, 141, 118, 116, 81, 33, 5, 1).
\]
\end{example}

\subsection{Normality}
Ein, Erman, and Lazarsfeld have conjectured that, for large values of $\bd$, the Betti numbers in any given row $\beta(\PP^1\times \PP^1,\bzero;\bd)$ should look approximately like a normal distribution~\cite[Conjecture~B]{EEL}.  Bruce proved that a similar phenomena holds for the first row when $\bd = (2,d_2)$ and $d_2\to \infty$ in~\cite[Theorem~A]{bruce-hirzebruch}, but that it fails for the second row under the same hypotheses~\cite[Theorem~B]{bruce-hirzebruch}.   See also~\cite{lemmensP1P1,erman-yang} for related results.

Our data, while somewhat limited, suggests that results similar to \cite[Theorem~A, Theorem~B]{bruce-hirzebruch} also hold for $\PP^1\times \PP^1$ embedded by $(3,d_2)$ as $d_2\to \infty$. In particular, as $d_{2}\to\infty$ the Betti numbers in the $q=1$ row of $\beta(\PP^1\times \PP^1,\bzero;(3,d_2))$ approach a normal distribution, while Betti numbers in the $q=2$ row do not. Figure~\ref{fig:3d normal distribution} highlights this for the $q=1$ row. 

It would be interesting to better understand what happens for the $q=2$ row and a fixed $d_{1}$. This is likely related to the phenomenon of asymptotic non-vanishing of syzygies in the semi-ample setting as discussed in \cite{bruce-semiample}.  Concretely, we ask:

\begin{question}
Does there exist $d_1\in \ZZ_{\geq2}$ such that the Betti numbers in the $q=2$ row of $\beta(\PP^1\times \PP^1,\bzero;(d_1,d_2))$ approach a normal distribution as $d_{2}\to \infty$? 
\end{question}

\section{Representation Theoretic Conjectures}\label{sec:rep theory conjectures}
Utilizing the representation theory of $\GL_2\times \GL_2$ provides the most concise way to express the syzygies of $\PP^1\times \PP^1$.  Our Schur functor data enabled us to make conjectures related to specific entries of the Betti tables. Additionally, our data raises questions regarding the ubiquity of redundant Schur functors.

\subsection{Specific Entries}
We first consider conjectures on specific $K_{p,q}$ groups.  
As noted earlier, the case when $\bb = \zero$ is of particular interest, as this case corresponds to the syzygies of the homogeneous coordinate ring of $\PP^1\times \PP^1$ under the embedding by $\cO_{\PP^{1}\times\PP^{1}}(\bd)$.  Moreover, based on our data and the unimodality conjectures from the previous section, we expect the extremal entries in a row to involve the fewest Schur functors.   

We thus are most interested in extremal entires in a row in the case $\bb=0$.  We first offer a conjecture about the last entry of the $q=1$ row:

\begin{conj}[Row $q=1$]\label{conj:last entry}
Let $\bd\in \ZZ^2_{\geq1}$ and $p=(d_1+1)(d_2-1)+d_1$.  (This is the largest value of $p$ such that $K_{p,1}(\PP^{1}\times\PP^{1},\bzero;\bd)\neq 0$ in this case.) Let 
\[
\ba \coloneqq\left(
\tbinom{d_1+1}{2}\tbinom{d_2}{1}\ , \ 
\tbinom{d_1+1}{2}\tbinom{d_2}{1}\ , \
\tbinom{d_{1}+1}{1}\tbinom{d_{2}+1}{2}-1\ , \
\tbinom{d_1+1}{1}\tbinom{d_2}{2}+1
\right)\in \ZZ^4.
\]
\begin{enumerate}
	\item  {\bf Last entry:} Assume $d_2>d_1$. Then $K_{p,1}(\PP^{1}\times\PP^{1},\bzero;\bd)$ is an irreducible Schur module.  Specifically, if $d_{2}>d_{1}$ then 
\begin{equation*}
K_{p,1}(\PP^{1}\times\PP^{1},\bzero;\bd)\cong \bS_{\ba+(0,0,-1,1)}.
\end{equation*}
	
\item  {\bf Second-to-last entry:}  Assume $d_2>d_1+1$.  
Then $K_{p-1,1}(\PP^{1}\times\PP^{1},\bzero;\bd)$ is the direct sum of $d_2$ distinct irreducible Schur modules.  Specifically, if $d_2>d_1+1$ then  
\[
K_{p-1,1}(\PP^{1}\times\PP^{1},\bzero;\bd)\cong \bigoplus_{i=0}^{d_2-1} \bS_{\ba+(0,-d_{1},-2-i,-d_{2}+2+i)}.
\]

\end{enumerate}
\end{conj}

Our next conjectures focus on the last entries in the $q=2$ row. In particular, the following conjecture describes the Schur functor decomposition for the last entry in the $q=2$ row for all $\bd$, as well as the decomposition for the second to last entry in the $q=2$ row in the special cases when $\bd=(2,d)$ and $\bd=(3,d)$.

\begin{conj}[Row $q=2$]\label{conj:q2-row}
Let $\bd\in \ZZ_{\geq1}^{2}$ and let $p=(d_{1}+1)(d_{2}+1)-3$. (This is the largest value of $p$ such that $K_{p,2}(\PP^{1}\times\PP^{1},\bzero;\bd)\neq0$.) 
\begin{enumerate}
	\item {\bf Last entry:}  The space $K_{p,2}(\PP^{1}\times\PP^{1},\bzero;\bd)$ is a unique irreducible Schur module.  Specifically, $K_{p,2}(\PP^{1}\times\PP^{1},\bzero;\bd)\cong\bS_{\ba}$, where
\begin{equation*}
\ba \coloneqq\left(
\tbinom{d_{1}+1}{2}\tbinom{d_{2}+1}{1}-1\ , \
\tbinom{d_{1}+1}{2}\tbinom{d_{2}+1}{1}-d_{1}+1 \ , \
\tbinom{d_{1}+1}{1}\tbinom{d_{2}+1}{2}-1 \ , \ 
\tbinom{d_{1}+1}{1}\tbinom{d_{2}+1}{2}-d_{2}+1 
\right)\in \ZZ^4.
\end{equation*}

\item {\bf Second-to-last entry, $\bd=(2,d)$:}  Assume that $\bd=(2,d)$. The space $K_{p,2}(\PP^{1}\times\PP^{1},\bzero;\bd)$ is the direct sum of  $d-2$ Schur modules.  Specifically, $K_{p-1,2}(\PP^{1}\times\PP^{1},\bzero;\bd)\cong \oplus_{i=0}^{d-3}\bS_{\ba+(0,0,-i,+i)}$, where 
\[
\ba \coloneqq\left(
3d+2 \ , \
3d \ , \
\tfrac{1}{2}(3d^{2}+3d-2)-1 \ , \ 
\tfrac{1}{2}(3d^{2}+3d-2) - 2(d_{2}-d_{1})-3
\right)\in \ZZ^4.
\]
\item {\bf Second-to-last entry, $\bd=(3,d)$:}  Assume that $\bd=(3,d)$. The space $K_{p,2}(\PP^{1}\times\PP^{1},\bzero;\bd)$ is the direct sum of $2d-3$  irreducible Schur module.  More specifically,  if

\begin{align*}
\ba &\coloneqq \left(
6d+5 \ , \ 
6d+1 \  , \ 
2d^{2}+2d-2 \ , \
2d^{2}+2d-2d+3
\right)\in \ZZ^4, \\
\bb &\coloneqq\left(
6d+4 \ , \ 
6d+2 \  , \ 
2d^{2}+2d-2 \ , \
2d^{2}+2d-2d+1
\right)\in \ZZ^4,
\end{align*}
then $K_{p-1,2}(\PP^{1}\times\PP^{1},\bzero;\bd)\cong \oplus_{i=0}^{d-3} \bS_{\ba+(0,0,-i,i)}\oplus \oplus_{j=0}^{d_{2}-2} \bS_{\bb+(0,0,-i,i)}$.
\end{enumerate}
\end{conj}

As we have only computed the full Betti table $\beta(\PP^1\times \PP^1,\bzero;(3,d_2))$ for four values of $d_2$, the evidence for part (3) of Conjecture~\ref{conj:q2-row} is admittedly scant. That said, the $\ba$'s in both parts (2) and (3) of Conjecture~\ref{conj:q2-row}, seem to fit into a potentially more general pattern. This leads us to ask the following question concerning the Schur functor decomposition for the second to last entry in the $q=2$ row in general. 

\begin{question}
Let $\bd\in \ZZ_{\geq1}^{2}$ and  let $p=(d_{1}+1)(d_{2}+1)-3$. (This is the largest value of $p$ such that $K_{p,2}(\PP^{1}\times\PP^{1},\bzero;\bd)\neq0$.) If 
\[
\ba \coloneqq \left( 
\tbinom{d_{1}+1}{2}\tbinom{d_{2}+1}{1}-1\ , \
\tbinom{d_{1}+1}{2}\tbinom{d_{2}+1}{1}-2d_{1}+1\ , \
\tbinom{d_{1}+1}{1}\tbinom{d_{2}+1}{2}-2\ , \
\tbinom{d_{1}+1}{1}\tbinom{d_{2}+1}{2}-2d_{2}+2
\right).
\]
then is it the case that as representations of $\GL_{2}\times \GL_{2}$:
\[
\bigoplus_{i=0}^{d_{2}-3} \bS_{\ba+(0,0,-i,i)} \subset K_{p-1,2}\left(\PP^{1}\times\PP^{1},\bzero;\bd\right)?
\]
\end{question}

\subsection{Redundant Schur Functors}
The central result of \cite{einLazarsfeld12} shows that asymptotically, Betti tables have numerous ``redundant'' entries.  That is, it is very often the case that both $K_{p,q}$ and $K_{p-1,q+1}$ will be nonzero.  These entries are ``redundant'' in the sense that they  could not be predicted by the Hilbert function of the module.  

A folklore question asks to find similar ``redundant'' representation in the Schur functor decomposition of $K_{p,q}$ and $K_{p-1,q+1}$.  More specifically, we consider examples of a Schur functor $\bS_{\mu}\otimes \bS_{\lambda}$ that appears in the Schur functor decomposition of both $K_{p,q}(\PP^1\times\PP^1, \bb;\bd)$ and  $K_{p-1,q+1}(\PP^1\times\PP^1, \bb;\bd)$.   In~\cite[Example 6.17 and Question 6.16]{bruceErmanGoldsteinYang18}, the authors give examples of redundant Schur functors for $\PP^{2}$ and ask whether redundant Schur functors occur frequently or sporadically.  Based upon our data, redundant Schur functors seem quite common for $\PP^1\times \PP^1$. For example, out of the approximately 200 pairs of $(\bb;\bd)$ for which we computed Schur functor computations rough two-thirds contained redundant Schur functors. 

While we did not find much of a pattern for when and where redundant Schur functors might occur, it would be interesting to explore that question further.  We did observe, anecdotally, that redundant Schur functors were more likely to occur if one of $d_1, d_2, b_1$ or $b_2$ is sufficiently large.  Focusing on the case when $\bb=\bzero$ our data suggests the following conjecture. 

\begin{conj}
If either $d_1$ or $d_2$ is sufficiently large, then there exists $p,q$ such that $K_{p,q}(\PP^1\times\PP^1,\bzero;\bd)$ has redundant Schur functors. 
\end{conj}

Furthermore, within each example, the number of redundant Schur functors seems able to be quite large both in terms of the total number and in terms of percentage of total Schur functors. For example, the largest total number of redundant Schur functors we observed is when $\bd=(2,10)$ and $\bb=(0,8)$; in this case, there are 596 redundant Schur functors out of 7135 total Schur functors (without multiplicity). The redundant Schur functors makes up the largest percentage of total Schur functors (counted without multiplicity) occurs when $\bd=(3,5)$ and $\bb=(2,4)$ where approximately $22.9\%$ of Schur functors are redundant. 

In addition, our data shows a number of examples where for a particular $p$ and $q$ all of the Schur functors appearing in the decomposition of $K_{p,q}(\PP^1\times\PP^1, \bb;\bd)$ are redundant. For example, when $\bd=(2,3)$ and $\bb=(1,2)$ both $K_{5,0}(\PP^1\times\PP^1, (1,2);(2,3))$ and $K_{4,1}(\PP^1\times\PP^1, (1,2);(2,3))$ are isomorphic to $\bS_{(8,3,11,6)}\oplus \bS_{(7,4,10,7)}\oplus\bS_{(6,5,9,8)}$ implying all of these Schur functors are redundant.  Appendix~\ref{apen:schur} includes the Schur functor decompositions of $K_{p,q}(\PP^1\times\PP^1, (1,2);(2,3))$ for all $p$ and $q$.

\section{Boij-S\"oderberg Theory Conjectures and Questions} \label{sec:boij sod}
\subsection{Background on Boij--S\"oderberg Theory}
Boij--S\"oderberg theory provides a way to decompose a Betti table as a positive rational sum of certain atomic building blocks called pure diagrams.  The theory was conjectured by~\cite{boij-sod} and the main results were proven in~\cite{eisenbud-schreyer-JAMS}.  See also \cite{floystad-expository,floystad-mccullough-peeva} for expository treatments of the theory or \cite{boij-sod2,eis-schrey2,eisenbud-erman,  floystad-multigraded, beks-local, beks-tensor,erman-semigroup, gibbons1, gibbons2} for more details on various aspects of the theory.

Having computed an array of Betti tables for embeddings of $\PP^1\times \PP^1$, we can analyze the pure diagrams and coefficients that arise in corresponding Boij--S\"oderberg decompositions. In order to get well-defined coefficients, we need to choose a specific set of representatives for the pure diagrams $\pi_{\bdelta}$.

Set $[n] = \{0,1,...,n-1\}$. 
Given a sequence of integers $\bdelta = (\delta_0,\ldots,\delta_r)$, called a \textit{degree sequence}, let $\pi_{\bdelta}$ be the Betti table with entries 
\begin{align*}
\beta_{i,j}(\pi_{\bdelta}) = \left\{
\begin{array}{cc}
\prod_{i\neq j} \frac{1}{|\delta_i-\delta_j|} & \text{ if } j= \delta_i \\
0 & \text{ if }j\neq \delta_i.
\end{array}
\right.
\end{align*}
For instance
\[
\pi_{(0,1,3,4)} = \begin{pmatrix}
\frac{1}{12} & \frac{1}{6} & \cdot&\cdot\\
\cdot&\cdot&\frac{1}{6}& \frac{1}{12}
\end{pmatrix}
\]
Note in particular, that $\pi_{\bdelta}$ will often have entries in $\QQ$, not in $\ZZ$.

For any graded Cohen-Macaulay module $M$ over a polynomial ring, there exists a unique set of degree sequences $C_M$ such that 
\begin{equation*}
\beta(M) = \sum_{\bdelta\in C_M} a_{\bdelta} \pi_{\bdelta} \hspace{20pt} \text{with } a_{\bdelta} \in \QQ.
\end{equation*}
This is called the \textit{Boij-S\"oderberg decomposition} of $M$, and the rational numbers $\{a_{\bdelta} \, | \, \bdelta\in C_M\}$ are called the \textit{Boij-S\"oderberg coefficients} of $M$.

\subsection{Conjectures on Boij-S\"oderberg coefficients}
Formulas for the coefficients have been found in certain cases where $M$ has a well-understood algebraic or combinatorial structure~\cite{ejo,mayes-tang,nagel-sturgeon,gibbons1,gibbons2,gibbons3,ees-filtering,erman-sam-supernatural}.  In this section, we aim to provide conjectures on Boij-S\"oderberg coefficients for the Betti tables of $\PP^1\times \PP^1$. 

One common feature of Boij-S\"oderberg decompositions, exhibited in many of the examples referenced above, is that they rarely ``skip over'' potential degree sequences. For instance, in the case $\bb=(0,0)$ and $\bd=(2,5)$, the shape of the Betti table is:
\[
\left(
\begin{array}{*{16}c}
*&\cdot&\cdot&\cdot&\cdot&\cdot&\cdot&\cdot&\cdot&\cdot&\cdot&\cdot&\cdot&\cdot&\cdot&\cdot\\
\cdot&*&*&*&*&*&*&*&*&*&*&*&*&*&*&\cdot\\
\cdot&\cdot&\cdot&\cdot&\cdot&\cdot&\cdot&\cdot&\cdot&\cdot&\cdot&\cdot&*&*&*&*
\end{array}
\right)
\]
where the zero entries are marked with $\cdot$ and nonzero entries are marked with $*$.
Based on this shape, there are only $4$ pure diagrams which could potentially arise in the Boij-S\"oderberg decomposition, depending on where you choose to shift from the 1st row to the 2nd row. (See also Example~\ref{ex:0025} below, which specifies the corresponding degree sequences.)  In this example, the coefficients of each such potential pure diagram turn out to be nonzero, although there is no obvious reason why this ought to be true.  Conjecture \ref{conj:BSSum} posits that this phenomenon occurs whenever $\bb=(0,0)$. More precisely, when $\bb=(0,0)$ and $\bd=(d_1,d_2)$ where $d_1\leq d_2$, the degree sequences which could possibly occur are given by
\begin{equation*}
\bdelta_j = 
[(d_1+1)(d_2+1)]\setminus \{1, (d_1+1)(d_2+1)-d_1-j\} \hspace{10pt} \text{ for } \hspace{10pt}  0\leq j\leq (d_1-1)(d_2-2)
\end{equation*}
and we conjecture the following.
\begin{conj}
	\label{conj:BSSum}
	The Boij-S\"oderberg coefficient $a_{\bdelta_j}$ is nonzero for each $j$.
\end{conj}

We now attempt to better understand the values of the nonzero coefficients.  Our first such conjecture, provides a complete description of the Boij-S\"oderberg coefficients in the case where $\bd=(2,d_2)$ and $\bb=(0,b_2)$ for $0\leq b_2 \leq d_2-2$ and $d_2\geq 3$.  In particular, taking $b_2=0$, this provides a complete conjectural description of the Boij-S\"oderberg coefficients of the homogeneous coordinate ring of $\PP^1\times \PP^1$ embedded by $\cO_{\PP^{1}\times\PP^{1}}(2,d_2)$.

\begin{conj}\label{conj:2d boij sod}
Let $\bd =(2,d_2)$ and $\bb=(0,b_2)$ for some $0\leq b_2 \leq d_2-2$.  Assume $d_2 \geq 3$.  The Boij-S\"oderberg decomposition will involve the degree sequences $\bdelta_j$ for $0\leq j \leq d_2-2$ where $\bdelta_j$ is defined as
	\begin{align*}
	\bdelta_j = \left\{ 
	\begin{array}{ll}
	[3(d_2+1)]\setminus \{b_2+1, 3d_2+1-j\} & 0\leq j\leq d_2-b_2-2  \\
	\left[3(d_2+1)\right]\setminus \{d_2-j-1, 2d_2+b_2+3\} & d_2-b_2-1 \leq j \leq d_2-2.
	\end{array} \right.
	\end{align*}
Moreover, the Boij-S\"oderberg coefficients corresponding to $\bdelta_j$ will be given by the formula
	\[a_{\bdelta_j}=\begin{cases}
	2(3d_2)! & j\neq d_2-b_2-2\\
	2(d_2+2)(3d_2)! & j=d_2-b_2-2.
	\end{cases}\]
In particular, all of the coefficients, except for the last one, will be identical, and as $d_2\to \infty$, the last coefficient will dominate.
\end{conj}
To prove Conjecture~\ref{conj:2d boij sod}, one might be able to use \cite[Corollary 5]{lemmens18}, which provides an explicit formula for the Betti numbers in this case.

\begin{example}\label{ex:0025}
As noted above, if we take $b_{2}=0$ and $d_2\geq 3$,  then Conjecture~\ref{conj:2d boij sod} implies that the Boij-S\"oderberg decomposition for the the homogeneous coordinate ring of $\PP^1\times \PP^1$ embedded by $\cO_{\PP^{1}\times\PP^{1}}(2,d_2)$ is:
\[
\beta\left(\PP^{1}\times \PP^{1},\bzero;(2,d_2)\right)=
2(3d_{2})!\left( \pi_{\bdelta_{0}}+\pi_{\bdelta_{1}}+\cdots + \pi_{\bdelta_{d_2-3}}\right) + 2(d_{2}+2)(3d_{2})! \pi_{\bdelta_{d_{2}-2}}
\]
where $\bdelta_{j}$ is the degree sequence $(0,2\ldots,3d_{2}-j,\widehat{3d_{2}+1-j} ,3d_{2}+2-j\ldots,3d_{2}+2)$. For example, if $d_{2}=5$ then we have:
\[
\begin{cases}
\bdelta_0 &=(0,2, 3, 4, 5, 6, 7, 8, 9, 10, 11, 12, 13, 14, 15, \widehat{16},17)\\
\bdelta_1 &=(0, 2, 3, 4, 5, 6, 7, 8, 9, 10, 11, 12, 13, 14, \widehat{15}, 16, 17)\\
\bdelta_2 &=(0, 2, 3, 4, 5, 6, 7, 8, 9, 10, 11, 12, 13, \widehat{14}, 15, 16, 17)\\
\bdelta_3 &=(0, 2, 3, 4, 5, 6, 7, 8, 9, 10, 11, 12, \widehat{13}, 14, 15, 16, 17)
\end{cases}
\]
and Conjecture~\ref{conj:2d boij sod} states that
\[
\beta\left(\PP^{1}\times \PP^{1},\bzero;(2,5)\right)= 2(15!)\left( \pi_{\bdelta_0} + \pi_{\bdelta_1}+\pi_{\bdelta_2}\right) +14(15!)\pi_{\bdelta_3}.
\]
\end{example}

\begin{remark}
Conjecture~\ref{conj:2d boij sod} would imply the following curious fact:  consider the Betti table of the homogeneous coordinate ring (that is, with $\bb=\zero$) of $\PP^1\times\PP^1$ embedded by $(2,d_2)$.  As $d_2\to \infty$, these Betti tables will be ``asymptotically pure'' in a sense that parallels the main result of~\cite{erman-high-degree}, where these Betti tables are asymptotically dominated by the contributions from a single pure diagram. See also~\cite{taylor,erman-sam-questions}.  It would be very interesting to better understand the limits under which such Betti tables are ``asymptotically pure''; this question is wide open for $\PP^2$ as well, as discussed in~\cite[\S6.3]{bruceErmanGoldsteinYang18}
\end{remark}

\begin{figure}
	\includegraphics[scale=0.6]{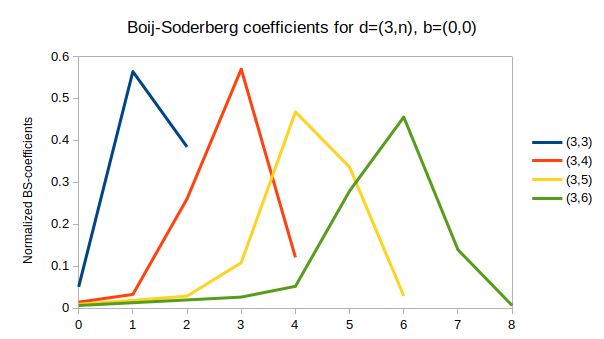}
	\caption{The Boij-Soderberg coefficients for $\bd=(3,n)$ and $\bb=(0,0)$} 
	\label{fig:BSb00d3n}
\end{figure}

When $\bd=(3,d_2)$ and $\bb=(0,0)$, we have a conjecture for roughly the first half of the coefficients.   Figure \ref{fig:BSb00d3n} displays these coefficients, rescaled by a factor of  $6d_2 (3d_2)!$ (so that these numbers sum to 1) to allow for a better comparison as $d_2$ grows. Notice that in each case, there is a set of small values followed by a peak.

\begin{conj}\label{conj:3d boij sod}
	For $\bb=(0,0)$ and $\bd=(3,d_2)$, with $d_2\geq 4$, the Boij-S\"oderberg coefficients for $j=0,\ldots, d_2-4$ are
	\begin{equation*}
	a_{\delta_j} = 
	\frac{(j+1)(4d_2+4)!}{4{4d_2+4 \choose 4}}.
	\end{equation*}
\end{conj}

\subsection{More questions}
Our data on Boij-S\"oderberg coefficients also illuminated some fascinating patterns which we were not able to convert into precise conjectures.  We conclude by drawing attention to a couple of these phenomena for curious readers.

As we saw in the previous conjectures, there are various situations where, if we fix some of the variables $b_1$, $b_2$, $d_1$, or $d_2$, then the number of Boij-S\"oderberg coefficients remains fixed.  When this happens, it is natural to understand how the individual coefficients depend on the remaining variables. 

The sum of the Boij-S\"oderberg coefficients of a module can be determined by the multiplicity of that module, and in the $\PP^1\times \PP^1$ case, this sum is
\begin{equation*}
\sum_{\bdelta\in C_{S(\bb;\bd)}}a_{\bdelta}=\frac{2d_1d_2}{((d_1+1)(d_2+1))^{\underline{3}}}\cdot ((d_1+1)(d_2+1))!
\end{equation*}
where $x^{\underline{n}}$ denotes the \textit{falling factorial}:
\begin{equation*}
x^{\underline{n}} = \prod_{k=0}^{n-1}(x+k).
\end{equation*}
To better analyze the coefficients, we rescale:
\begin{equation*}
b_{\bdelta} = \frac{a_{\bdelta}}{((d_1+1)(d_2+1))!} \text{ and note that } \sum_{\bdelta\in C_{S(\bb;\bd)}}b_{\bdelta}=\frac{2d_1d_2}{((d_1+1)(d_2+1))^{\underline{3}}}.
\end{equation*}
While it appears difficult to give concrete conjectural formulae for the Boij-S\"oderberg coefficients for larger values of $\bd$ than those studied in the previous subsection, the above equation suggests something about the behavior of the $b_{\bdelta}$ as a rational function of $d_1$ and/or of $d_2$.

The following is a concrete conjecture in this direction:

\begin{conj}\label{conj:bs-rational}
	For $\bb=(d_1-1,b_2)$ and $\bd=(d_1,d_2)$, with $d_1\leq d_2$, $0\leq b_2\leq d_2-2$, the degree sequences appearing in the Boij-S\"oderberg decomposition of $\beta(\PP^1\times \PP^1, \bb; \bd)$ are precisely:
	\begin{equation*}
	\delta_j = 
	[(d_1+1)(d_2+1)-1]\setminus \{(b_2+1)d_1-j\} \hspace{20pt}  0\leq j\leq (b_2+1)(d_1-1).
	\end{equation*}
	For any fixed $b_2$ and $d_1$, the coefficient $a_{\delta_j}$ has the form
	\begin{equation*}
	a_{\delta_j}=p_j(d_2)((d_1+1)(d_2+1))!
	\end{equation*}
	where $p_{j}$ is some degree $-2$ rational function in $d_2$.
\end{conj}
Evidence for this conjecture is provided in Appendix~\ref{appendix:BS}.  Moreover, the above discussion and Conjecture~\ref{conj:bs-rational} suggest the following question.

\begin{question}
What is the value of 
\[
\lim_{d_2\rightarrow \infty}a_{\delta_j}\frac{((d_1+1)(d_2+1))^{\underline{3}}}{2d_1d_2((d_1+1)(d_2+1))!}?
\]
\end{question}

\begin{table}[]
	\centering
	\begin{tabular}{|c|c|l|}
		\hline
		$d_1$ & $\bb$ & Normalized BS-coefficients as $d_2\to \infty$  \\
		\hline
		\hline
		$2$  & $(1,0)$  & $(1/2,1/2)$ \\
		& $(1,1)$  & $(2/9,5/9,2/9)$ \\
		& $(1,2)$  & $(8/81,32/81,67/162,5/54)$ \\
		\hline
		$3$  & $(2,0)$  & $(3/8,1/4,3/8)$ \\
		\hline
		$4$  & $(3,0)$  & $(8/25,4/25,11/50,3/10)$ \\
		\hline
	\end{tabular}
	\caption{Asymptotic values of BS-coefficients as $d_2\to \infty$, normalized so that the numbers sum to 1}
	\label{tab:asymptoticValues}
\end{table}

Most of the conjectures considered in this section can be understood as being motivated by the following overarching but vague question:

\begin{question}
To what extent, and under what additional restrictions, can the Boij-S\"oderberg coefficients of $\beta(\PP^1\times \PP^1,\bb;\bd)$ be understood as rational functions in $b_1$, $b_2$, $d_1$ and/or $d_2$?
\end{question}

We end with a mystery.  In Appendix~\ref{appendix:BS}, we plotted the Boij-S\"oderberg coefficients of  $\beta(\PP^1\times \PP^1,\bb;\bd)$, after rescaling so that the sum of the coefficients is $1$, for various natural families depending on $\bb$ and/or $\bd$.  We simply note that the graphics suggest a remarkable uniformity among these families as one varies the parameters.  Can one explain, or even precisely describe, this phenomenon?

\newpage

\appendix

\section{Number and size of matrices computed}
\label{appendix:matrices}

We record, for nearly all $(\bb;\bd)$ pairs for which we have complete data, the number of matrices in the relevant range and the size of the largest matrix. 
{\footnotesize 
	\begin{table}[tbh!]
		\centering
		\begin{tabular}{|l||l|l|l||l|l|l|} 
			\hline
			$\bd$ & $\bb$ &  Number of& Largest & $\bb$ &  Number of& Largest  \\
			& & matrices & matrix & & matrices & matrix\\
			\hline \hline
			$(2,4)$ &	$(0,0)$ & $75$ & $625 \times 2431$ & $(1,0)$ &0 & N/A \\
			&	$(0,1)$ & 0 & N/A & $(1,1)$ & $23$ & $73\times 81$ \\
			&	$(0,2)$ & $17$ & $19\times 15 $ & $(1,2)$ &$109$ & $554\times 909$ \\
			&	$(0,3)$ & $44$ & $77\times 82$ & $(1,3)$ &$212$ & $1387\times 3171$ \\
			\hline
			$(2,5)$&	$(0,0)$ & $216$ & $3386 \times 13946$ & $(1,0)$ &0 & N/A \\
			&$(0,1)$ &$101$ & $1508 \times 6988$ & $(1,1)$ &$31$ &$108\times 116$ \\
			&$(0,2)$ &$20$ & $23\times 18$ & $(1,2)$ &$135$ & $1245\times 1911$ \\
			&$(0,3)$ &$55$ & $116 \times 117$ & $(1,3)$ &$297$ & $5302\times 10822$ \\
			&$(0,4)$ &$110$ & $434 \times 552$ &$(1,4)$ &$486$ & $9432\times 25262$ \\
			\hline
			$(2,6)$&$(0,0)$ & $466$ & $18,902\times 81,386$ & $(1,0)$ & 0& N/A  \\
			&$(0,1)$ &$273$ & $8547\times 40,922$ & $(1,1)$ &$35$ &$148\times 156$ \\
			&$(0,2)$ &$150$ & $3075\times 16,649$ & $(1,2)$ &$171$ & $2476\times 3607$ \\
			&$(0,3)$ &$62$ & $159\times 155$ &$(1,3)$ &$367$ & $15,588\times 29,403$ \\
			&$(0,4)$ &$131$ & $723\times 868$ &$(1,4)$ &$651$ & $44,886\times 107,138$ \\
			&$(0,5)$ &$212$ & $2512\times 3580$ &$(1,5)$ &$919$ & $62,250\times 187,699$ \\
			\hline
			$(2,7)$&$(0,0)$ & $831$ & $108,060\times 482,053$ &$(1,0)$ &0 &NA  \\
			&$(0,1)$ & $573$ & $49,808\times 243,840$ & $(1,1)$ &$43$ &$196\times 204$ \\
			&$(0,2)$ & $368$ & $18,682\times 102,154$ & $(1,2)$ &$197$ & $4392\times 6171$ \\
			&$(0,3)$ & $226$ & $5600\times 34,800$ & $(1,3)$ &$456$ & $39,140\times 69,452$ \\
			&$(0,4)$ & $148$ & $1118\times 1286$ &$(1,4)$ &$795$ & $163,325\times 358,383$ \\
			&$(0,5)$ & $251$ & $4562\times 6132$ &$(1,5)$ &$1198$ & $352,746\times 949,098$ \\
			&$(0,6)$ & $385$ & $14,782\times 22,836$ & $(1,6)$ &$1619$ & $436,912\times 1,248,208$ \\
			\hline
			$(2,8)$&$(0,0)$ & $1391$ & $627,537\times 2,886,389$ & $(0,7)$ & $622$ & $87,266\times 144,514$ \\
			&$(0,1)$ & $995$ & $291,943\times 1,460,756$ & $(1,0)$ & 0  & NA  \\
			&$(0,2)$ & $721$ & $113,886\times 627,766$ & $(1,1)$ & $47$ & $249\times 255$  \\
			&$(0,3)$ & $479$ & $36,350\times 224,623$ & $(1,2)$ & $233$ & $7310\times 9966$  \\
			&$(0,4)$ & $348$ & $9408\times 66,110$ & $(1,3)$ & $527$ & $86,245\times 146,042$  \\
			&$(0,5)$ & $280$ & $7594\times 9764$ & $(1,4)$ & $968$ & $498,024\times 1,022,361$ \\
			&$(0,6)$ & $445$ & $28,470\times 41,648$ &  & &  \\
			\hline
			$(2,9)$ & $(0,7)$ & $809$ & $177,658 \times 278,759$ & $(1,2)$ & $259$ & $11,415\times 15,196$ \\
			&$(1,0)$ & $0$ & N/A & $(1,3)$ &  $616$  &  $174,144 \times 283,727$  \\
			&$(1,1)$ & $55$ & $310 \times 314 $ &  &  &   \\ 
			\hline
			$(2,10)$ & $(0,8)$ & $1106$ & $1,111,726\times 1,843,366$ & $(1,2)$ & $295$ & $17,132 \times 22,350$ \\
			&$(1,0)$ & $0$ & N/A & $(1,3)$ & $687$ & $325,114 \times 513,364$   \\
			&$(1,1)$ & $59$ & $376 \times 378$ &  &  &   \\ 
			\hline
			$(2,11)$ &$(1,0)$ & $0$ & N/A  & $(1,2)$ & $321$ & $24,649 \times 31,638$ \\
			&$(1,1)$ & $67$ & $ 450\times 450 $ & $(1,3)$ & $776$ & $574,112\times 882,626$   \\
			\hline
		\end{tabular}
		\caption{Matrix data}
		\label{table:matrixData}
	\end{table}

	\begin{table}
		\centering
		\begin{tabular}{|l||l|l|l||l|l|l|} 
			\hline
			$\bd$ & $\bb$ &  Number of& Largest & $\bb$ &  Number of& Largest  \\
			& & matrices & matrix & & matrices & matrix\\
			\hline \hline
			$(3,3)$ & $(0,0)$ & $104$ & $1772\times 6180$ & $(1,1)$ & $31$ & $88\times 96$  \\
			&$(0,1)$ & $0$ & NA & 	$(1,2)$ & $125$ & $740\times 1204$ \\
			&$(0,2)$ & $19$ & $20\times 16$ & $(2,2)$ & $308$ & $2838\times 7308$ \\
			\hline
			$(3,4)$ &$(0,0)$ & $521$ & $25,320\times 87,114$ & $(1,2)$ & $177$ & $2038\times 3023$ \\
			&$(0,1)$ & $148$ & $6678\times 29,840$ & $(1,3)$ & $368$ & $11,086\times 21,945$ \\
			&$(0,2)$ & $23$ & $26\times 20$ & $(2,0)$ & $24$ & $24\times 20$ \\
			&$(1,3)$ & $58$ & $130\times 140$ & $(2,1)$ & $164$ & $1956\times 2975$ \\	 
			&$(1,0)$ & $122$ & $2407\times 12,740$  & $(2,2)$ & $468$ & $19,478\times 43,618$ \\
			&$(1,1)$ & $37$ & $138\times 144$ & $(2,3)$ & $836$ & $35,556\times 96,730$ \\
			\hline
			$(3,5)$&$(0,0)$ & $1344$ & $361,276\times 1,231,276$ & $(1,3)$ & $495$ & $39,424\times 70,894$ \\
			&$(0,1)$ & $711$ & $119,254\times 505,443$ & $(1,4)$ & $858$ & $162,286\times 367,093$ \\
			&$(0,2)$ & $236$ & $19168\times 104246$ & $(2,0)$ & $29$ & $30\times 24$ \\
			&$(0,3)$ & $72$ & $196\times 200$ & $(2,1)$ & $219$ & $4350\times 6320$ \\
			&$(0,4)$ & $153$ & $1078\times 1328$ & $(2,2)$ & $618$ & $87,401\times 178,536$ \\
			&$(1,0)$ & $338$ & $18,014\times 101,895$ &  $(2,3)$ & $1217$ & $348,702\times 971,100$ \\
			&$(1,1)$ & $210$ & $5758\times 35,668$ & $(2,4)$ & $1891$ & $467,124\times 1,322,104$ \\
			&$(1,2)$ & $215$ & $4470\times 6248$ & & & \\
			\hline 
			$(3,6)$&$(0,3)$ & $334$ & $45,094\times 290,746$ & $(1,3)$ & $598$ & $110,702\times 186,050$ \\
			&$(0,4)$ & $181$ & $1774\times 2076$ & $(1,4)$ & $1106$ & $697,950\times 1,436,165$ \\
			&$(0,5)$ & $296$ & $8224\times 11,390$ & $(2,0)$ & $34$ & $34\times 28$ \\
			&$(1,0)$ & $740$ & $142,906\times 845,408$ & $(2,1)$ & $259$ & $8347\times 11,760$ \\	 
			&$(1,1)$ & $476$ & $44,876\times 290,369$ & $(2,2)$ & $793$ & $300,091\times 573,890$ \\
			&$(1,2)$ & $468$ & $11,665\times 83,466$ & & & \\
			\hline
			$(3,7)$ & $(0,4)$ & $518$ & $94,088\times 700,128$ & $(1,2)$ & $827$ & $97,064\times 709,416$ \\
			&$(0,5)$ & $348$ & $14,768 \times 19,378$ & $(1,3)$ & $968$ &  $268940 \times 428636$  \\
			&$(0,6)$ & $554$ & $68,616 \times 102,344$ & $(2,0)$ & $39$  & $40\times 32$   \\ 
			&$(1,0)$ & $1130$ & $1,128,854\times 6,980,468$ & $(2,1)$ & $314$  & $14,766\times 20,308$   \\ 
			&$(1,1)$ & $965$ & $369,576 \times 2,450,184$ & $(2,2)$ & $950$  & $855,136\times 1,556,128$   \\ 
			\hline
			$(3,8)$ & $(0,5)$ & $745$ & $177,432\times 1,500,926$ & $(2,0)$ & $44$ & $44\times 36$ \\
			&$(0,6)$ & $636$ & $130,144\times 184,592$ & $(2,1)$ & $354$ & $24,136\times 32,574$  \\
			&$(0,7)$ & $900$ & $ 553,291 \times 879,321 $ & $(2,2)$ & $1130$ &  $2,124,896\times 3,719,448$ \\ 
			\hline
			$(3,9)$ & $(0,7)$ & $1026$ & $1,105,918\times 1,673,092$ & $(2,1)$ & $409$ & $37,620\times 49,992$ \\
			&$(2,0)$ & $49$ & $50\times 40$ & $ $ &  &   \\
			\hline
						$(4,4)$ & $(0,0)$ & $1715$ & $853,068\times 2,722,820$ & $(1,2)$ & $228$ & $5269\times 7364$ \\
			&$(0,1)$ & $764$ & $165,929\times 743,227$ &	$(1,3)$ & $501$ & $50,156\times 91,458$  \\
			&$(0,2)$ & $198$ & $6518\times 43,768$ & $(2,2)$ & $682$ & $121,747\times 241,924$\\
			&$(0,3)$ & $80$ & $207\times 222$ & $(2,3)$ & $1321$ & $581,410\times 1,582,730$  \\
			& $(1,1)$ & $249$ & $24,765\times 138,553$ & & & \\
			\hline
			$(4,5)$ & $(0,2)$ & $1102 $ & $731,824\times 4,033,789 $ & $(2,0)$ & $525 $ & $63,634\times 456,031 $ \\
			&$(0,3)$ & $315 $ & $14,409\times 117,520 $ &	$(2,1)$ & $521 $ & $15,511\times 120,826 $  \\
			&$(0,4)$ & $207 $ & $2121\times 2596 $ & $(2,2)$ & $879 $ & $517,511\times 948,896 $\\
			&$(1,1)$ & $1075 $ & $755,881\times 4,074,383 $ & $(2,3)$ & $ $ & $ $  \\
			&$(1,2)$ & $559 $ & $70,246\times 471,986 $ &$(3,0)$ &$94 $ & $300\times 318$\\
			&$(1,3)$ & $661 $ & $171,904\times 287,389$ &$(3,1)$ &$673 $ & $169,940\times 292,984$ \\
			\hline
			$(4,6)$ & $(0,4)$ & $507$ & $27,864\times 267,592$ & $(2,2)$ & $1381$ & $1713790\times 2964636$  \\
			& $(0,5)$ & $411$ & $21,318\times 28,941$ & $(3,0)$ & $116$ & $417\times 438$ \\
			& $(1,3)$ & $1132$ & $471,259\times 740,692$ & $(3,1)$ & $813$ & $462,729\times 767,366$ \\
			& $(2,1)$ & $956$ & $157,164\times 1,277,412$ &  &  &  \\
			\hline		
			$(4,7)$ & $(0,5)$ & $791$ & $49,046\times 542,194$ & $(3,0)$ & $130$ & $540\times 568$  \\
			& $(0,6)$ & $762$ & $221,972\times 324,448$ & $(3,1)$ & $991$ & $1,100,334\times 1,771,080$ \\
			\hline
		\end{tabular}
		\caption{Matrix data}
		\label{table:matrixData3}
	\end{table}

\mbox{~}

\begin{landscape}

\section{Total Betti Numbers}

{\fontsize{3.5}{4.2}
\setcounter{MaxMatrixCols}{20}

\begin{flalign*}
&\beta(\bzero;(2,2))=
\begin{array}{*{30}c}
&0&1&2&3&4&5&6&7\\ 
 	\text{0:}&1&\text{.}&\text{.}&\text{.}&\text{.}&\text{.}&\text{.}&\text{.}\\ 
 	\text{1:}&\text{.}&20&64&90&64&20&\text{.}&\text{.}\\
 	\text{2:}&\text{.}&\text{.}&\text{.}&\text{.}&\text{.}&\text{.}&1&\text{.}\\
 	\end{array}&
\end{flalign*}

\begin{flalign*}
&\beta(\bzero;(2,3))=
\begin{array}{*{30}c}
&0&1&2&3&4&5&6&7&8&9&10\\ 
 	\text{0:}&1&\text{.}&\text{.}&\text{.}&\text{.}&\text{.}&\text{.}&\text{.}&\text{.}&\text{.}&\text{.}\\
 	\text{1:}&\text{.}&43&222&558&840&798&468&147&8&\text{.}&\text{.}\\ 
 	\text{2:}&\text{.}&\text{.}&\text{.}&\text{.}&\text{.}&\text{.}&\text{.}&\text{.}&9&2&\text{.}\\
\end{array}&
\end{flalign*}

\begin{flalign*}
&\beta(\bzero;(2,4))=
\begin{array}{*{30}c}
&0&1&2&3&4&5&6&7&8&9&10&11&12&13\\
 	\text{0:}&1&\text{.}&\text{.}&\text{.}&\text{.}&\text{.}&\text{.}&\text{.}&\text{.}&\text{.}&\text{.}&\text{.}&\text{.}&\text{.}\\ 
 	\text{1:}&\text{.}&75&536&1947&4488&7095&7920&6237&3344&1089&120&11&\text{.}&\text{.}\\
 	\text{2:}&\text{.}&\text{.}&\text{.}&\text{.}&\text{.}&\text{.}&\text{.}&\text{.}&\text{.}&\text{.}&66&24&3&\text{.}\\
 	\end{array}&
\end{flalign*}

\begin{flalign*}
&\beta(\bzero;(2,5))=
\begin{array}{*{30}c}
&0&1&2&3&4&5&6&7&8&9&10&11&12&13&14&15&16\\
\text{0:}&1&\text{.}&\text{.}&\text{.}&\text{.}&\text{.}&\text{.}&\text{.}&\text{.}&\text{.}&\text{.}&\text{.}&\text{.}&\text{.}&\text{.}&\text{.}&\text{.}\\
\text{1:}&\text{.}&116&1060&5040&15652&34580&56628&70070&65780&46332&23660&8008&1260&195&14&\text{.}&\text{.}\\
\text{2:}&\text{.}&\text{.}&\text{.}&\text{.}&\text{.}&\text{.}&\text{.}&\text{.}&\text{.}&\text{.}&\text{.}&\text{.}&455&210&45&4&\text{.}\\
\end{array}&
\end{flalign*}

\begin{flalign*}
&\beta(\bzero;(2,6))=
\begin{array}{*{30}c}
&0&1&2&3&4&5&6&7&8&9&10&11&12&13&14&15&16&17&18&19\\ 
 	\text{0:}&1&\text{.}&\text{.}&\text{.}&\text{.}&\text{.}&\text{.}&\text{.}&\text{.}&\text{.}&\text{.}&\text{.}&\text{.}&\text{.}&\text{.}&\text{.}&\text{.}&\text{.}&\text{.}&\text{.}\\ 
 	\text{1:}&\text{.}&166&1848&10863&42432&120360&259488&436254&579904&612612&512720&335478&166464&58344&11424&2295&288&17&\text{.}&\text{.}\\
 	\text{2:}&\text{.}&\text{.}&\text{.}&\text{.}&\text{.}&\text{.}&\text{.}&\text{.}&\text{.}&\text{.}&\text{.}&\text{.}&\text{.}&\text{.}&3060&1632&459&72&5&\text{.}\\
 	\end{array}&
\end{flalign*}

\begin{flalign*}
&\beta(\bzero;(2,7))=
\begin{array}{*{30}c}
&0&1&2&3&4&5&6&7&8&9&10&11&12&13&14&15&16&17&18&19&20&21&22\\ 
 	\text{0:}&1&\text{.}&\text{.}&\text{.}&\text{.}&\text{.}&\text{.}&\text{.}&\text{.}&\text{.}&\text{.}&\text{.}&\text{.}&\text{.}&\text{.}&\text{.}&\text{.}&\text{.}&\text{.}&\text{.}&\text{.}&\text{.}&\text{.}\\ 
 	\text{1:}&\text{.}&225&2954&20685&97356&337155&901170&1912806&3281680&4598874&5290740&4996810&3852744&2393430&1166676&421515&95760&22610&3780&399&20&\text{.}&\text{.}\\ 
 	\text{2:}&\text{.}&\text{.}&\text{.}&\text{.}&\text{.}&\text{.}&\text{.}&\text{.}&\text{.}&\text{.}&\text{.}&\text{.}&\text{.}&\text{.}&\text{.}&\text{.}&20349&11970&3990&840&105&6&\text{.}\\
 	\end{array}&
\end{flalign*}

\begin{flalign*}
&\beta(\bzero;(2,8))=
\begin{array}{*{30}c}
&0&1&2&3&4&5&6&7&8&9&10&11&12&13&14&15&16&17&18&19&20&21&22&23&24&25\\ 
 	\text{0:}&1&\text{.}&\text{.}&\text{.}&\text{.}&\text{.}&\text{.}&\text{.}&\text{.}&\text{.}&\text{.}&\text{.}&\text{.}&\text{.}&\text{.}&\text{.}&\text{.}&\text{.}&\text{.}&\text{.}&\text{.}&\text{.}&\text{.}&\text{.}&\text{.}&\text{.}\\ 
 	\text{1:}&\text{.}&293&4432&36018&198352&811118&2586672&6628853&13921072&24270543&35421472&43474508&44930592&39017108&28289632&16915833&8152672&3023603&765072&201894&40480&5796&528&23&\text{.}&\text{.}\\ 
 	\text{2:}&\text{.}&\text{.}&\text{.}&\text{.}&\text{.}&\text{.}&\text{.}&\text{.}&\text{.}&\text{.}&\text{.}&\text{.}&\text{.}&\text{.}&\text{.}&\text{.}&\text{.}&\text{.}&134596&85008&31878&8096&1380&144&7&\text{.}\\
 	\end{array}&
\end{flalign*}

\begin{flalign*}&
\beta(\bzero;(3,3))=
\begin{array}{*{30}c}
&0&1&2&3&4&5&6&7&8&9&10&11&12&13&14\\ 
 	\text{0:}&1&\text{.}&\text{.}&\text{.}&\text{.}&\text{.}&\text{.}&\text{.}&\text{.}&\text{.}&\text{.}&\text{.}&\text{.}&\text{.}&\text{.}\\ 
 	\text{1:}&\text{.}&87&676&2691&6864&12155&15444&14157&9152&3861&780&22&\text{.}&\text{.}&\text{.}\\ 
 	\text{2:}&\text{.}&\text{.}&\text{.}&\text{.}&\text{.}&\text{.}&\text{.}&\text{.}&\text{.}&\text{.}&165&144&39&4&\text{.}\\
 	\end{array}&
\end{flalign*}

\begin{flalign*}
&\beta(\bzero;(3,4))=
\begin{array}{*{30}c}
&0&1&2&3&4&5&6&7&8&9&10&11&12&13&14&15&16&17&18\\ 
 	\text{0:}&1&\text{.}&\text{.}&\text{.}&\text{.}&\text{.}&\text{.}&\text{.}&\text{.}&\text{.}&\text{.}&\text{.}&\text{.}&\text{.}&\text{.}&\text{.}&\text{.}&\text{.}&\text{.}\\ 
 	\text{1:}&\text{.}&147&1530&8364&30192&78540&153816&232050&272272&247962&172380&87516&28560&3939&238&15&\text{.}&\text{.}&\text{.}\\ 
 	\text{2:}&\text{.}&\text{.}&\text{.}&\text{.}&\text{.}&\text{.}&\text{.}&\text{.}&\text{.}&\text{.}&\text{.}&\text{.}&1287&3094&1800&528&85&6&\text{.}\\
 	\end{array}&
\end{flalign*}

\begin{flalign*}
&\beta(\bzero;(3,5))=
\begin{array}{*{30}c}
&0&1&2&3&4&5&6&7&8&9&10&11&12&13&14&15&16&17&18&19&20&21&22\\ 
 	\text{0:}&1&\text{.}&\text{.}&\text{.}&\text{.}&\text{.}&\text{.}&\text{.}&\text{.}&\text{.}&\text{.}&\text{.}&\text{.}&\text{.}&\text{.}&\text{.}&\text{.}&\text{.}&\text{.}&\text{.}&\text{.}&\text{.}&\text{.}\\ 
 	\text{1:}&\text{.}&223&2912&20265&94696&325185&860472&1804278&3049120&4191894&4702880&4291378&3147312&1805570&759696&195390&25088&3247&360&19&\text{.}&\text{.}&\text{.}\\ 
 	\text{2:}&\text{.}&\text{.}&\text{.}&\text{.}&\text{.}&\text{.}&\text{.}&\text{.}&\text{.}&\text{.}&\text{.}&\text{.}&\text{.}&\text{.}&6435&37856&41684&20520&6270&1240&147&8&\text{.}\\
 	\end{array}
&\end{flalign*}

\begin{flalign*}
&\beta((1,1);(2,3))=
\begin{array}{*{45}c}
&0&1&2&3&4&5&6&7&8&9&10&11\\
\text{0:}&4&28&72&56&\text{.}&\text{.}&\text{.}&\text{.}&\text{.}&\text{.}&\text{.}&\text{.}\\
\text{1:}&\text{.}&\text{.}&8&168&504&672&528&252&68&8&\text{.}&\text{.}\\
\text{2:}&\text{.}&\text{.}&\text{.}&\text{.}&\text{.}&\text{.}&\text{.}&\text{.}&\text{.}&\text{.}&\text{.}&\text{.}\\
\end{array}&
\end{flalign*}

\begin{flalign*}
&\beta((1,1);(2,4))=
\begin{array}{*{45}c}
&0&1&2&3&4&5&6&7&8&9&10&11&12&13&14\\
\text{0:}&4&36&120&120&\text{.}&\text{.}&\text{.}&\text{.}&\text{.}&\text{.}&\text{.}&\text{.}&\text{.}&\text{.}&\text{.}\\
\text{1:}&\text{.}&\text{.}&32&660&2772&5808&7920&7524&5060&2376&744&140&12&\text{.}&\text{.}\\
\text{2:}&\text{.}&\text{.}&\text{.}&\text{.}&\text{.}&\text{.}&\text{.}&\text{.}&\text{.}&\text{.}&\text{.}&\text{.}&\text{.}&\text{.}&\text{.}\\
\end{array}&
\end{flalign*}

\begin{flalign*}
&
\beta((1,1);(2,5))=
\begin{array}{*{45}c}
&0&1&2&3&4&5&6&7&8&9&10&11&12&13&14&15&16&17\\
\text{0:}&4&44&180&220&\text{.}&\text{.}&\text{.}&\text{.}&\text{.}&\text{.}&\text{.}&\text{.}&\text{.}&\text{.}&\text{.}&\text{.}&\text{.}&\text{.}\\
\text{1:}&\text{.}&\text{.}&80&1820&9828&28028&54340&77220&82940&68068&42588&20020&6860&1620&236&16&\text{.}&\text{.}\\
\text{2:}&\text{.}&\text{.}&\text{.}&\text{.}&\text{.}&\text{.}&\text{.}&\text{.}&\text{.}&\text{.}&\text{.}&\text{.}&\text{.}&\text{.}&\text{.}&\text{.}&\text{.}&\text{.}\\
\end{array}
&
\end{flalign*}

\begin{flalign*}
&
\beta((1,1);(2,6))=
\begin{array}{*{45}c}
&0&1&2&3&4&5&6&7&8&9&10&11&12&13&14&15&16&17&18&19&20\\
\text{0:}&4&52&252&364&\text{.}&\text{.}&\text{.}&\text{.}&\text{.}&\text{.}&\text{.}&\text{.}&\text{.}&\text{.}&\text{.}&\text{.}&\text{.}&\text{.}&\text{.}&\text{.}&\text{.}\\
\text{1:}&\text{.}&\text{.}&160&4080&26928&97104&243984&461448&680680&797368&747864&562224&337008&159120&57936&15708&2988&356&20&\text{.}&\text{.}\\
\text{2:}&\text{.}&\text{.}&\text{.}&\text{.}&\text{.}&\text{.}&\text{.}&\text{.}&\text{.}&\text{.}&\text{.}&\text{.}&\text{.}&\text{.}&\text{.}&\text{.}&\text{.}&\text{.}&\text{.}&\text{.}&\text{.}\\
\end{array}&
\end{flalign*}

}

{\fontsize{2.5}{3}
\begin{flalign*}
&
\beta((1,1);(2,7))=
\begin{array}{*{45}c}
&0&1&2&3&4&5&6&7&8&9&10&11&12&13&14&15&16&17&18&19&20&21&22&23\\
\text{0:}&4&60&336&560&\text{.}&\text{.}&\text{.}&\text{.}&\text{.}&\text{.}&\text{.}&\text{.}&\text{.}&\text{.}&\text{.}&\text{.}&\text{.}&\text{.}&\text{.}&\text{.}&\text{.}&\text{.}&\text{.}&\text{.}\\
\text{1:}&\text{.}&\text{.}&280&7980&62244&271320&837216&1976760&3708040&5643456&7054320&7289464&6240360&4418640&2573664&1220940&464436&138320&31080&4956&500&24&\text{.}&\text{.}\\
\text{2:}&\text{.}&\text{.}&\text{.}&\text{.}&\text{.}&\text{.}&\text{.}&\text{.}&\text{.}&\text{.}&\text{.}&\text{.}&\text{.}&\text{.}&\text{.}&\text{.}&\text{.}&\text{.}&\text{.}&\text{.}&\text{.}&\text{.}&\text{.}&\text{.}\\
\end{array}&
\end{flalign*}

\begin{flalign*}
&\beta((1,1);(2,8))=
\begin{array}{*{45}c}
&0&1&2&3&4&5&6&7&8&9&10&11&12&13&14&15&16&17&18&19&20&21&22&23&24&25&26\\
\text{0:}&4&68&432&816&\text{.}&\text{.}&\text{.}&\text{.}&\text{.}&\text{.}&\text{.}&\text{.}&\text{.}&\text{.}&\text{.}&\text{.}&\text{.}&\text{.}&\text{.}&\text{.}&\text{.}&\text{.}&\text{.}&\text{.}&\text{.}&\text{.}&\text{.}\\
\text{1:}&\text{.}&\text{.}&448&14168&127512&651728&2384272&6749028&15363172&28765088&44930592&59075408&65731792&62047008&49685152&33668228&19208772&9152528&3598672&1147608&289432&55568&7632&668&28&\text{.}&\text{.}\\
\text{2:}&\text{.}&\text{.}&\text{.}&\text{.}&\text{.}&\text{.}&\text{.}&\text{.}&\text{.}&\text{.}&\text{.}&\text{.}&\text{.}&\text{.}&\text{.}&\text{.}&\text{.}&\text{.}&\text{.}&\text{.}&\text{.}&\text{.}&\text{.}&\text{.}&\text{.}&\text{.}&\text{.}\\
\end{array}&
\end{flalign*}

\begin{flalign*}
&\beta((1,1);(2,9))=
\begin{array}{*{45}c}
&0&1&2&3&4&5&6&7&8&9&10&11&12&13&14&15&16&17&18&19&20&21&22&23&24&25&26&27&28&29\\
\text{0:}&4&76&540&1140&\text{.}&\text{.}&\text{.}&\text{.}&\text{.}&\text{.}&\text{.}&\text{.}&\text{.}&\text{.}&\text{.}&\text{.}&\text{.}&\text{.}&\text{.}&\text{.}&\text{.}&\text{.}&\text{.}&\text{.}&\text{.}&\text{.}&\text{.}&\text{.}&\text{.}&\text{.}\\
\text{1:}&\text{.}&\text{.}&672&23400&238680&1399320&5920200&19536660&52295100&116233260&217809540&347677200&476050320&561632400&572330160&504131940&383467500&251213820&141098100&67490280&27232920&9149400&2513160&549900&92196&11124&860&32&\text{.}&\text{.}\\
\text{2:}&\text{.}&\text{.}&\text{.}&\text{.}&\text{.}&\text{.}&\text{.}&\text{.}&\text{.}&\text{.}&\text{.}&\text{.}&\text{.}&\text{.}&\text{.}&\text{.}&\text{.}&\text{.}&\text{.}&\text{.}&\text{.}&\text{.}&\text{.}&\text{.}&\text{.}&\text{.}&\text{.}&\text{.}&\text{.}&\text{.}\\
\end{array}&
\end{flalign*}
}

{\fontsize{1.1}{1.32}
\begin{flalign*}
&\beta((1,1);(2,10))=
\begin{array}{*{45}c}
&0&1&2&3&4&5&6&7&8&9&10&11&12&13&14&15&16&17&18&19&20&21&22&23&24&25&26&27&28&29&30&31&32\\
\text{0:}&4&84&660&1540&\text{.}&\text{.}&\text{.}&\text{.}&\text{.}&\text{.}&\text{.}&\text{.}&\text{.}&\text{.}&\text{.}&\text{.}&\text{.}&\text{.}&\text{.}&\text{.}&\text{.}&\text{.}&\text{.}&\text{.}&\text{.}&\text{.}&\text{.}&\text{.}&\text{.}&\text{.}&\text{.}&\text{.}&\text{.}\\
\text{1:}&\text{.}&\text{.}&960&36540&416556&2755116&13232700&49877100&153476700&394877340&863111340&1620609900&2634716700&3729663900&4614746220&5002540020&4756176900&3965381700&2895246900&1846402740&1024391940&491645700&202562100&70913700&20805876&5020596&970340&144420&15540&1076&36&\text{.}&\text{.}\\
\text{2:}&\text{.}&\text{.}&\text{.}&\text{.}&\text{.}&\text{.}&\text{.}&\text{.}&\text{.}&\text{.}&\text{.}&\text{.}&\text{.}&\text{.}&\text{.}&\text{.}&\text{.}&\text{.}&\text{.}&\text{.}&\text{.}&\text{.}&\text{.}&\text{.}&\text{.}&\text{.}&\text{.}&\text{.}&\text{.}&\text{.}&\text{.}&\text{.}&\text{.}\\
\end{array}&
\end{flalign*}

\begin{flalign*}
&\beta((1,1);(2,11))=
\begin{array}{*{45}c}
&0&1&2&3&4&5&6&7&8&9&10&11&12&13&14&15&16&17&18&19&20&21&22&23&24&25&26&27&28&29&30&31&32&33&34&35\\
\text{0:}&4&92&792&2024&\text{.}&\text{.}&\text{.}&\text{.}&\text{.}&\text{.}&\text{.}&\text{.}&\text{.}&\text{.}&\text{.}&\text{.}&\text{.}&\text{.}&\text{.}&\text{.}&\text{.}&\text{.}&\text{.}&\text{.}&\text{.}&\text{.}&\text{.}&\text{.}&\text{.}&\text{.}&\text{.}&\text{.}&\text{.}&\text{.}&\text{.}&\text{.}\\
\text{1:}&\text{.}&\text{.}&1320&54560&687456&5063168&27214528&115345296&401097840&1172439840&2928294720&6322199520&11900027040&19651420800&28603734720&36819120360&42004911960&42523491120&38211096000&30459702240&21507388320&13418545920&7371224640&3548173200&1487147376&538278048&166451648&43353376&9329760&1614976&216128&20988&1316&40&\text{.}&\text{.}\\
\text{2:}&\text{.}&\text{.}&\text{.}&\text{.}&\text{.}&\text{.}&\text{.}&\text{.}&\text{.}&\text{.}&\text{.}&\text{.}&\text{.}&\text{.}&\text{.}&\text{.}&\text{.}&\text{.}&\text{.}&\text{.}&\text{.}&\text{.}&\text{.}&\text{.}&\text{.}&\text{.}&\text{.}&\text{.}&\text{.}&\text{.}&\text{.}&\text{.}&\text{.}&\text{.}&\text{.}&\text{.}\\
\end{array}&
\end{flalign*}

\begin{flalign*}
&\beta((1,1);(2,12))=
\begin{array}{*{45}c}
&0&1&2&3&4&5&6&7&8&9&10&11&12&13&14&15&16&17&18&19&20&21&22&23&24&25&26&27&28&29&30&31&32&33&34&35&36&37&38\\
\text{0:}&4&100&936&2600&\text{.}&\text{.}&\text{.}&\text{.}&\text{.}&\text{.}&\text{.}&\text{.}&\text{.}&\text{.}&\text{.}&\text{.}&\text{.}&\text{.}&\text{.}&\text{.}&\text{.}&\text{.}&\text{.}&\text{.}&\text{.}&\text{.}&\text{.}&\text{.}&\text{.}&\text{.}&\text{.}&\text{.}&\text{.}&\text{.}&\text{.}&\text{.}&\text{.}&\text{.}&\text{.}\\
\text{1:}&\text{.}&\text{.}&1760&78540&1083852&8796480&52312128&246256560&954881840&3125556896&8781000480&21428722224&45830660400&86489553600&144765466560&215756224200&287156386440&341989309200&364915966800&349058361960&299274762600&229802523840&157793918400&96668031600&52670597616&25418685600&10807648544&4021262960&1298064240&359506752&84194880&16352028&2563260&311640&27576&1580&44&\text{.}&\text{.}\\
\text{2:}&\text{.}&\text{.}&\text{.}&\text{.}&\text{.}&\text{.}&\text{.}&\text{.}&\text{.}&\text{.}&\text{.}&\text{.}&\text{.}&\text{.}&\text{.}&\text{.}&\text{.}&\text{.}&\text{.}&\text{.}&\text{.}&\text{.}&\text{.}&\text{.}&\text{.}&\text{.}&\text{.}&\text{.}&\text{.}&\text{.}&\text{.}&\text{.}&\text{.}&\text{.}&\text{.}&\text{.}&\text{.}&\text{.}&\text{.}\\
\end{array}&
\end{flalign*}
}

{\fontsize{2}{2.4}
\begin{flalign*}
&\beta((1,1);(3,3))=
\begin{array}{*{45}c}
&0&1&2&3&4&5&6&7&8&9&10&11&12&13&14\\
\text{0:}&4&39&144&165&\text{.}&\text{.}&\text{.}&\text{.}&\text{.}&\text{.}&\text{.}&\text{.}&\text{.}&\text{.}&\text{.}\\
\text{1:}&\text{.}&\text{.}&22&780&3861&9152&14157&15444&12155&6864&2691&676&87&\text{.}&\text{.}\\
\text{2:}&\text{.}&\text{.}&\text{.}&\text{.}&\text{.}&\text{.}&\text{.}&\text{.}&\text{.}&\text{.}&\text{.}&\text{.}&\text{.}&1&\text{.}\\
\end{array}&
\end{flalign*}

\begin{flalign*}
&\beta((1,1);(3,4))=
\begin{array}{*{45}c}
&0&1&2&3&4&5&6&7&8&9&10&11&12&13&14&15&16&17&18\\
\text{0:}&4&50&240&364&\text{.}&\text{.}&\text{.}&\text{.}&\text{.}&\text{.}&\text{.}&\text{.}&\text{.}&\text{.}&\text{.}&\text{.}&\text{.}&\text{.}&\text{.}\\
\text{1:}&\text{.}&\text{.}&58&2448&16728&57120&132600&228072&301444&311168&251940&159120&77112&27744&6936&1020&30&\text{.}&\text{.}\\
\text{2:}&\text{.}&\text{.}&\text{.}&\text{.}&\text{.}&\text{.}&\text{.}&\text{.}&\text{.}&\text{.}&\text{.}&\text{.}&\text{.}&\text{.}&\text{.}&\text{.}&16&2&\text{.}\\
\end{array}&
\end{flalign*}

\begin{flalign*}
&\beta((1,1);(3,5))=
\begin{array}{*{45}c}
&0&1&2&3&4&5&6&7&8&9&10&11&12&13&14&15&16&17&18&19&20&21&22\\
\text{0:}&4&61&360&680&\text{.}&\text{.}&\text{.}&\text{.}&\text{.}&\text{.}&\text{.}&\text{.}&\text{.}&\text{.}&\text{.}&\text{.}&\text{.}&\text{.}&\text{.}&\text{.}&\text{.}&\text{.}&\text{.}\\
\text{1:}&\text{.}&\text{.}&127&6020&52269&233016&721905&1697688&3155710&4739056&5819814&5878600&4888282&3333360&1846914&817836&281295&71288&11795&756&39&\text{.}&\text{.}\\
\text{2:}&\text{.}&\text{.}&\text{.}&\text{.}&\text{.}&\text{.}&\text{.}&\text{.}&\text{.}&\text{.}&\text{.}&\text{.}&\text{.}&\text{.}&\text{.}&\text{.}&\text{.}&\text{.}&\text{.}&190&40&3&\text{.}\\
\end{array}&
\end{flalign*}

\begin{flalign*}
&\beta((1,1);(3,6))=
\begin{array}{*{45}c}
&0&1&2&3&4&5&6&7&8&9&10&11&12&13&14&15&16&17&18&19&20&21&22&23&24&25&26\\
\text{0:}&4&72&504&1140&\text{.}&\text{.}&\text{.}&\text{.}&\text{.}&\text{.}&\text{.}&\text{.}&\text{.}&\text{.}&\text{.}&\text{.}&\text{.}&\text{.}&\text{.}&\text{.}&\text{.}&\text{.}&\text{.}&\text{.}&\text{.}&\text{.}&\text{.}\\
\text{1:}&\text{.}&\text{.}&240&12600&132480&728640&2823480&8424900&20189400&39801960&65523780&90930960&106977600&106977600&90930960&65523780&39801960&20189400&8424900&2823480&728640&132480&12600&1124&48&\text{.}&\text{.}\\
\text{2:}&\text{.}&\text{.}&\text{.}&\text{.}&\text{.}&\text{.}&\text{.}&\text{.}&\text{.}&\text{.}&\text{.}&\text{.}&\text{.}&\text{.}&\text{.}&\text{.}&\text{.}&\text{.}&\text{.}&\text{.}&\text{.}&\text{.}&2024&552&72&4&\text{.}\\
\end{array}&
\end{flalign*}

\begin{flalign*}
&\beta((1,1);(3,7))=
\begin{array}{*{45}c}
&0&1&2&3&4&5&6&7&8&9&10&11&12&13&14&15&16&17&18&19&20&21&22&23&24&25&26&27&28&29&30\\
\text{0:}&4&83&672&1771&\text{.}&\text{.}&\text{.}&\text{.}&\text{.}&\text{.}&\text{.}&\text{.}&\text{.}&\text{.}&\text{.}&\text{.}&\text{.}&\text{.}&\text{.}&\text{.}&\text{.}&\text{.}&\text{.}&\text{.}&\text{.}&\text{.}&\text{.}&\text{.}&\text{.}&\text{.}&\text{.}\\
\text{1:}&\text{.}&\text{.}&408&23548&290493&1900080&8838765&31962060&93776865&228914400&472526145&833976780&1268123745&1669794480&1909884465&1900189620&1644131655&1234857120&802110855&447943860&213176535&85322640&28144935&7396740&1448811&175392&20678&1564&57&\text{.}&\text{.}\\
\text{2:}&\text{.}&\text{.}&\text{.}&\text{.}&\text{.}&\text{.}&\text{.}&\text{.}&\text{.}&\text{.}&\text{.}&\text{.}&\text{.}&\text{.}&\text{.}&\text{.}&\text{.}&\text{.}&\text{.}&\text{.}&\text{.}&\text{.}&\text{.}&\text{.}&\text{.}&20475&6552&1134&112&5&\text{.}\\
\end{array}&
\end{flalign*}

\begin{flalign*}
&\beta((1,1);(4,4))=
\begin{array}{*{45}c}
&0&1&2&3&4&5&6&7&8&9&10&11&12&13&14&15&16&17&18&19&20&21&22&23\\
\text{0:}&4&64&400&816&\text{.}&\text{.}&\text{.}&\text{.}&\text{.}&\text{.}&\text{.}&\text{.}&\text{.}&\text{.}&\text{.}&\text{.}&\text{.}&\text{.}&\text{.}&\text{.}&\text{.}&\text{.}&\text{.}&\text{.}\\
\text{1:}&\text{.}&\text{.}&112&6776&64064&304304&1003200&2515524&5002624&8072416&10708096&11757200&10708096&8072416&5002624&2515524&1003200&304304&64064&6776&112&\text{.}&\text{.}&\text{.}\\
\text{2:}&\text{.}&\text{.}&\text{.}&\text{.}&\text{.}&\text{.}&\text{.}&\text{.}&\text{.}&\text{.}&\text{.}&\text{.}&\text{.}&\text{.}&\text{.}&\text{.}&\text{.}&\text{.}&\text{.}&816&400&64&4&\text{.}\\
\end{array}&
\end{flalign*}

}

\end{landscape}

\newpage

\section{Example of Schur Functor Decomposition}\label{apen:schur}

\begin{align*}
K_{0,0}(\PP^{1}\times\PP^{1},(1,2);(2,3))\;=\;&\bS_{(1,0,2,0)}\\
\\
K_{1,0}(\PP^{1}\times\PP^{1},(1,2);(2,3))\;=\;&\bS_{(2,1,3,2)}\oplus\bS_{(2,1,4,1)}\oplus\bS_{(2,1,5,0)}\oplus\bS_{(3,0,3,2)}\oplus\bS_{(3,0,4,1)}\\
\\
K_{2,0}(\PP^{1}\times\PP^{1},(1,2);(2,3))\;=\;&\bS_{(3,2,4,4)}\oplus\bS_{(3,2,5,3)}^{\oplus2}\oplus\bS_{(3,2,6,2)}^{\oplus2}\oplus\bS_{(3,2,7,1)}\oplus\bS_{(4,1,4,4)}\oplus\bS_{(4,1,5,3)}^{\oplus2}\\
&\oplus\bS_{(4,1,6,2)}^{\oplus2}\oplus\bS_{(4,1,7,1)}\oplus\bS_{(5,0,5,3)}\\
\\
K_{3,0}(\PP^{1}\times\PP^{1},(1,2);(2,3))\;=\;&\bS_{(4,3,6,5)}^{\oplus2}\oplus\bS_{(4,3,7,4)}^{\oplus3}\oplus\bS_{(4,3,8,3)}^{\oplus2}\oplus\bS_{(4,3,9,2)}\oplus\bS_{(5,2,6,5)}^{\oplus2}\oplus\bS_{(5,2,7,4)}^{\oplus3}\\
&\oplus\bS_{(5,2,8,3)}^{\oplus2}\oplus\bS_{(5,2,9,2)}\oplus\bS_{(6,1,6,5)}\oplus\bS_{(6,1,7,4)}\oplus\bS_{(6,1,8,3)}\\
\\
K_{4,0}(\PP^{1}\times\PP^{1},(1,2);(2,3))\;=\;&\bS_{(5,4,7,7)}\oplus\bS_{(5,4,8,6)}^{\oplus2}\oplus\bS_{(5,4,9,5)}^{\oplus2}\oplus\bS_{(5,4,10,4)}\oplus\bS_{(5,4,11,3)}\oplus\bS_{(6,3,7,7)}
\\
&\oplus\bS_{(6,3,8,6)}^{\oplus2}\oplus\bS_{(6,3,9,5)}^{\oplus2}\oplus\bS_{(6,3,10,4)}\oplus\bS_{(7,2,8,6)}\oplus\bS_{(7,2,9,5)}\oplus\bS_{(7,2,10,4)}\\
\\
K_{4,1}(\PP^{1}\times\PP^{1},(1,2);(2,3))\;=\;&\bS_{(6,5,9,8)}\oplus\bS_{(7,4,10,7)}\oplus\bS_{(8,3,11,6)}\\
\\
K_{5,0}(\PP^{1}\times\PP^{1},(1,2);(2,3))\;=\;&\bS_{(6,5,9,8)}\oplus\bS_{(7,4,10,7)}\oplus\bS_{(8,3,11,6)}\\
\\
K_{5,1}(\PP^{1}\times\PP^{1},(1,2);(2,3))\;=\;&\bS_{(7,6,10,10)}\oplus\bS_{(7,6,11,9)}^{\oplus2}\oplus\bS_{(7,6,12,8)}^{\oplus2}\oplus\bS_{(7,6,13,7)}\oplus\bS_{(7,6,14,6)}\oplus\bS_{(8,5,10,10)}\\
&\oplus\bS_{(8,5,11,9)}^{\oplus2}\oplus\bS_{(8,5,12,8)}^{\oplus2}\oplus\bS_{(8,5,13,7)}\oplus\bS_{(9,4,11,9)}\oplus\bS_{(9,4,12,8)}\oplus\bS_{(9,4,13,7)}\\
\\
K_{6,1}(\PP^{1}\times\PP^{1},(1,2);(2,3))\;=\;&\bS_{(8,7,12,11)}^{\oplus2}\oplus\bS_{(8,7,13,10)}^{\oplus3}\oplus\bS_{(8,7,14,9)}^{\oplus2}\oplus\bS_{(8,7,15,8)}\oplus\bS_{(9,6,12,11)}^{\oplus2}\oplus\bS_{(9,6,13,10)}^{\oplus3}\\
&\oplus\bS_{(9,6,14,9)}^{\oplus2}\oplus\bS_{(9,6,15,8)}\oplus\bS_{(10,5,12,11)}\oplus\bS_{(10,5,13,10)}\oplus\bS_{(10,5,14,9)}\\
\\
K_{7,1}(\PP^{1}\times\PP^{1},(1,2);(2,3))\;=\;&\bS_{(9,8,13,13)}\oplus\bS_{(9,8,14,12)}^{\oplus2}\oplus\bS_{(9,8,15,11)}^{\oplus2}\oplus\bS_{(9,8,16,10)}\oplus\bS_{(10,7,13,13)}\\
&\oplus\bS_{(10,7,14,12)}^{\oplus2}\oplus\bS_{(10,7,15,11)}^{\oplus2}\oplus\bS_{(10,7,16,10)}\oplus\bS_{(11,6,14,12)}\\
\\
K_{8,1}(\PP^{1}\times\PP^{1},(1,2);(2,3))\;=\;&\bS_{(10,9,15,14)}\oplus\bS_{(10,9,16,13)}\oplus\bS_{(10,9,17,12)}\oplus\bS_{(11,8,15,14)}\oplus\bS_{(11,8,16,13)}\\
\\
K_{9,1}(\PP^{1}\times\PP^{1},(1,2);(2,3))\;=\;&\bS_{(11,10,17,15)}
\end{align*}

\newpage

\section{Boij-S\"oderberg coefficients for $\bd = (d_1,d_2)$, $\bb=(d_1-1,b_2)$}
\label{appendix:BS}

We record the  Boij-S\"oderberg coefficients for $\bd = (d_1,d_2)$, $\bb=(d_1-1,b_2)$, normalized so that the coefficients sum to 1. This provides evidence for Conjecture \ref{conj:bs-rational} and illustrates the asymptotic behavior of the Boij-S\"oderberg coefficients in 1-parameter families of fixed degree sequence length.

\begin{multicols}{2}
	\includegraphics[scale=0.4]{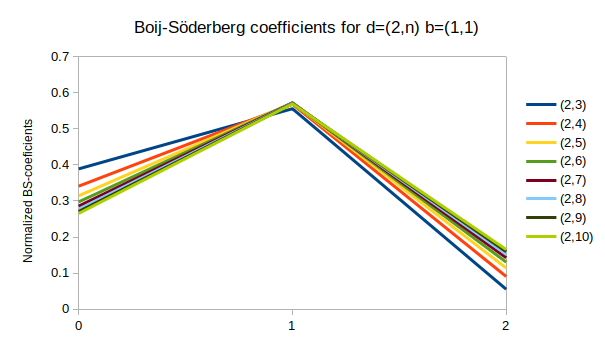}\\
	\includegraphics[scale=0.4]{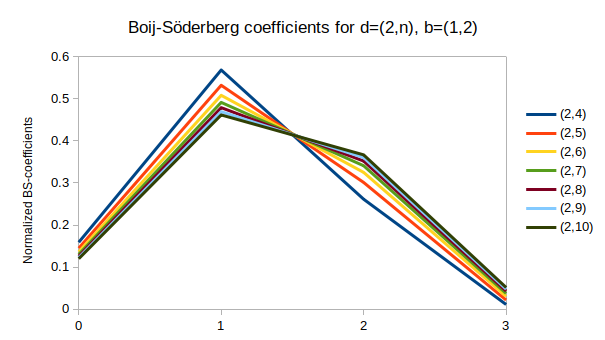}\\
	\includegraphics[scale=0.4]{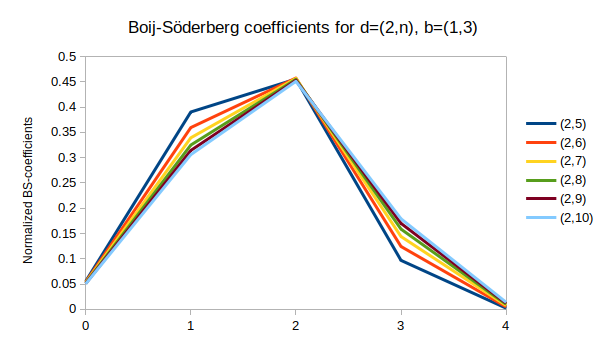}\\
	\includegraphics[scale=0.4]{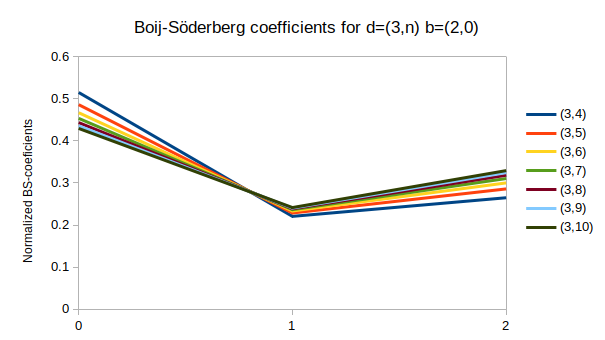}\\
	\includegraphics[scale=0.4]{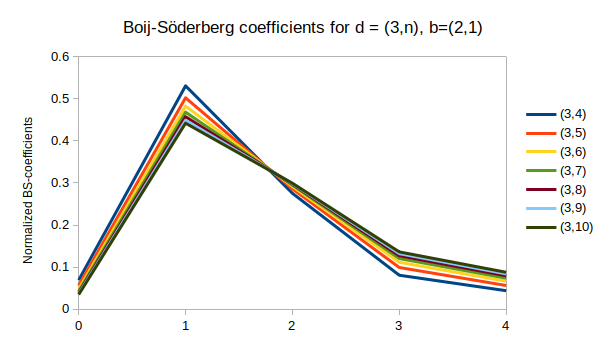}\\
	\includegraphics[scale=0.4]{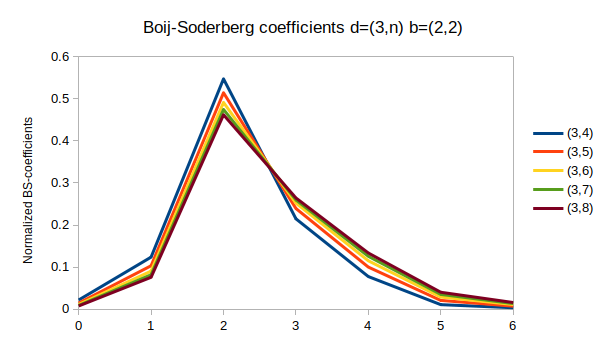}\\
	\includegraphics[scale=0.4]{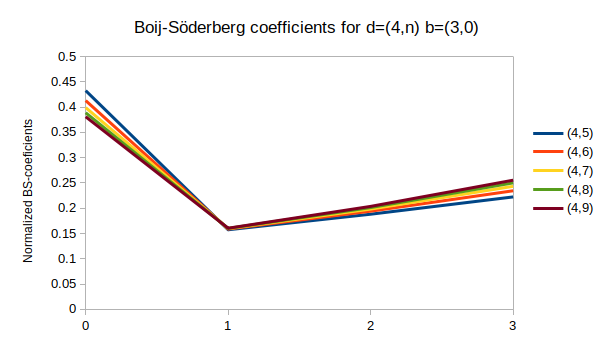}\\
	\includegraphics[scale=0.4]{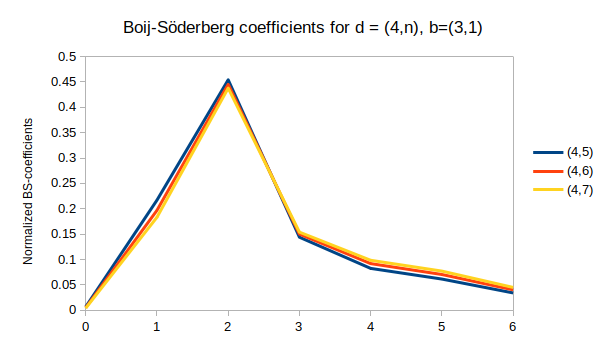}
\end{multicols}

\end{document}